\documentstyle[aps]{revtex}
\begin{document}

\title{\bf Extensions and Results from a Method for
Evaluating Fractional Integrals}



\author{\bf Victor Kowalenko}

\address{School of Physics\\
University of Melbourne, Parkville, Victoria 3052,
Australia,}

\vspace{2.0cm}

\author{\bf and M.L. Glasser}

\address{Department of Physics and Department of
Mathematics and Computer Science\\
Clarkson University, Potsdam, NY 13699-5815, USA.}


\maketitle
\begin{abstract}
We present a method derived from Laplace transform theory that 
enables the evaluation of fractional integrals. This method is
adapted and extended in a variety of ways to demonstrate  
its utility in deriving alternative representations
for other classes of integrals. We also use the method in conjunction
with several different techniques to derive many results that have not
appeared in tables of integrals.
\end{abstract}

\section{\bf Introduction}

In studying a heat diffusion problem recently one of the 
authors (M.L.G.) encountered the following fractional integral 
\begin{equation}
I= \int_1^{\infty} dx \;{e^{a^2x^2}erfc(ax) \over
\sqrt{x^2-1}} \quad. 
\label{one}\end{equation}
This integral does not appear in a succinct form in the more familiar 
tables of integrals such as Gradshteyn and Ryzhik \cite{gra80} 
and Prudnikov et al \cite{pru86} although the latter gives
a more general version of it in terms of $_2F_2$ and $_1F_1$
hypergeometric functions (see No. 2.8.7.6). So, a new method
was devised for evaluating it, unbeknownst to the authors at the
time that the succinct form they eventually obtained 
was given on p. 139 of Apelblat's book \cite{ape83}. Since this simple
method turns out to be very effective for evaluating fractional
integrals, we aim here to present the technique with many new results
to demonstrate its utility. In addition, we aim to extend the method
in different directions to illustrate its utility
and in so doing, we shall present many results not published
previously. 

As we shall see, our technique can be applied to
fractional integrals of the form
\begin{equation}
I_{\mu}(p)= \int_p^{\infty} dx\; {F(x) \over (x-p)^{\mu}} \quad,
\label{two}\end{equation}
with the aid of tables of LTs such as those by
Oberhettinger and Badii \cite{obe73}, and by Prudnikov et al 
\cite{pru92}-\cite{pru92a} and
extended to many other classes of integrals, 
not of a fractional form. The technique can also be used
in the numerical computation of integrals, since a slowly 
convergent integral can be transformed into a rapidly convergent
one.

\section{\bf The Method}

\vskip .1in
\noindent
Theorem: Given that the inverse LT for $F(x)$ in Eq.\ (\ref{two}) 
exists and is denoted by $G(t)$, then the integral $I_{\mu}(p)$
for $p$ real and greater than zero can be written as
\begin{equation}
I_{\mu}(p)= \Gamma(1-\mu) \int_0^{\infty}dt\;t^{\mu-1}
e^{-pt}G(t) \quad,
\label{three}\end{equation}
provided ${\rm Re}\; \mu < 1$ and $G(t)$ is $O(t^{\alpha})$ as $t \to 0+$
where ${\rm Re}\; \alpha > -{\rm Re}\; \mu$ while as $t \to \infty$ $G(t)$
is $O(t^{\beta}e^{pt})$ where ${\rm Re}\; \beta < -{\rm Re}\; \mu$.

Proof: Since the inverse LT for $F(x)$ exists, Eq.\ (\ref{two}) can
be rewritten as
\begin{equation}
I_{\mu}(p)= \int_p^{\infty} dx \;{1 \over (x-p)^{\mu}}\;
\int_0^{\infty}dt\; \exp(-xt) G(t) \quad.
\label{four}\end{equation}
Taking note of the conditions on $G(t)$, we can interchange
the order of the integrations whereupon we note that the
$x$-integral is simply the LT of $(x-p)^{-\mu} \Theta(x-p)$,
where $\Theta(x)$ represents the Heaviside step-function. Now from
p. 21 of Ref.\ \cite{obe73} we have
\begin{equation}
{\cal L}\bigl\{(x-p)^{-\mu}\Theta(x-p)\bigr\}= \int_p^{\infty}
dx\; e^{-xt}(x-p)^{-\mu}=
t^{\mu-1}e^{-pt}\Gamma(1-\mu)\quad ,
\label{five}\end{equation}
where ${\rm Re}\; \mu < 1$ and ${\rm Re}\; x > 0$.[ The condition on $\mu$
was overlooked on p. 17 of Prudnikov et al \cite{pru92}.]
Introducing Eq.\ (\ref{five}) into Eq.\ (\ref{four}) then yields
Eq.\ (\ref{three}), i.e. $I_{\mu}(p)= \Gamma(1-\mu)
{\cal L}_p\bigl\{t^{\mu-1}G(t)\bigr\}$, where the subscript $p$
denotes the variable with which the LT is taken.

Before developing further results from Eq.\ (\ref{three}) we now
consider some examples. Let us first consider the known
fractional integral 
\begin{equation}
I_{\mu}(p)=\int_p^{\infty}dx\;{x^{\alpha}(x+a)^{-\nu} \over
(x-p)^{\mu}} \quad,
\label{six}\end{equation}
which appears as No. 2.2.6.24 in \cite{pru86a}.
Thus, $F(x) = x^{\alpha}(x+a)^{-\nu}$ and its inverse LT is given
in Ref.\ \cite{pru92a} as
\begin{equation}
{\cal L}^{-1}\{x^{\alpha}(x+a)^{-\nu}\}={t^{-\alpha+\nu-1} \over
\Gamma(\nu-\alpha)} \; _1F_1(\nu;-\alpha+\nu;-at) \quad,
\label{seven}\end{equation}
where ${\rm Re} (\alpha-\nu) < 0$, ${\rm Re}\; x > \{0, -{\rm Re}\; a\}$
for $\alpha, -\nu \neq 0,1,2,...,\; {\rm Re}\; x > -{\rm Re}\;a$ for
$\alpha =0,1,2,...$ and ${\rm Re}\; x > 0$ for $-\nu = 0,1,2,...$. 
Then using Eq.\ (\ref{three}) we find that
\begin{equation}
I_{\mu}(p)= {\Gamma(1-\mu) \over \Gamma(\nu-\alpha)} 
\int_0^{\infty}dt\; t^{\mu+\nu-\alpha-2} e^{-pt}\;
_1F_1(\nu;\nu-\alpha;-at) \quad .
\label{eight}\end{equation}
The above integral is given on p. 510 of Ref.\ \cite{pru92}
and hence, Eq.\ (\ref{eight}) becomes
\begin{equation}
I_{\mu}(p)= {\Gamma(1-\mu) \Gamma(\mu+\nu-\alpha-1) \over
\Gamma(\nu-\alpha)\; p^{\mu+\nu-\alpha-1}}\;
_2F_1(\nu,\mu+\nu-\alpha-1;\nu-\alpha;-a/p)\;\;=\end{equation}
\begin{equation}
I_{\mu}(p)= (p+a)^{-\nu} p^{1-\mu+\alpha}\;B(1-\mu,\mu+\nu-\alpha-1)
\;_2F_1(1-\mu, \nu; \nu-\alpha; (1+p/a)^{-1}) \quad,
\label{onezero}\end{equation}
where $0 < {\rm Re} (1-\mu) < {\rm Re} (\nu-\mu)$. 

For the special situation where $\alpha=0$, we find that Eq.(9)
reduces to
\begin{equation}
I_{\mu}(p)={\Gamma(1-\mu) \over \Gamma(\nu)}\; {\Gamma(\mu + \nu-1)
\over (p+a)^{\mu+\nu-1}} \quad,
\label{oneone}\end{equation}
where we have utilised No. 7.3.1.1 from Ref.\ \cite{pru90}. For
$\mu=1/2$ and $p=1$, the above result reduces to No. 3.196.2
in Gradshteyn and Ryzhik \cite{gra80}. 
For the special case where $\mu =0$ we find that $I_{\mu}(p)=
\beta^{-\nu}\gamma(\nu,\beta/(a+p))$, where $\gamma(\nu,x)$ is
the incomplete gamma function, while if $\mu=1-\nu/2$ and ${\rm Re}\;\nu<2$,
then
\begin{equation}
I_{\mu}(p)={\sqrt{\pi}\; \Gamma(\nu/2)\; \beta^{-(\nu-1)/2}\;
e^{-\beta/2(a+p)} \over 2 (a+p)^{3/2}}\;\Bigl[
I_{(\nu-1)/2}\Bigl({\beta \over 2(a+p)}\Bigr)-
I_{(\nu+1)/2}\Bigl({\beta \over 2(a+p)}\Bigr)\Bigr] \quad.
\label{onefive}\end{equation}

The above analysis has demonstrated that our theorem can
be used during the intermediate steps of an evaluation to
obtain important results such as Eqs.(9) and (11).
To obtain 
results  not given the tables of integrals, first put $\nu=2n+\nu,\;
\mu=-\mu-n$ in (11), then multiply both 
sides of the original integral by $(-1)^n \beta^n/n!$.
After interchanging summations and integrations, we get
\begin{eqnarray}&
\int_0^{\infty}dt\;{t^{\mu} e^{-\beta t/(t+p)^2} \over (t+p)^{\nu}}
={\Gamma(\mu+1) \Gamma(\nu-\mu-1) \over p^{\nu-\mu-1} \Gamma(\nu)}\;
_2F_2\left(\mu+1,\nu-\mu-1;{\nu \over 2},{\nu+1 \over 2};-
{\beta \over 4 p}\right)\;\;,
\label{oneeightc}&\end{eqnarray}
where ${\rm Re}\;\mu>-1$.  This is equivalent to the interesting
formula 
$$\int_0^1x^{\mu}(1-x)^{\nu-\mu-2}e^{-\beta
x(1-x)}dx=\frac{\Gamma(\mu+1)
\Gamma(\nu-\mu-1)}{\Gamma(\nu)}\:_2F_2(\mu+1,\nu-\mu-1;\frac{\nu}{2},
\frac{\nu+1}{2};-\beta).$$
Similarly, we also find
\begin{eqnarray}&
\int_0^{\infty}dt\;{t^{\mu-\alpha/2} \over (t+p)^{\nu-\alpha}}\;
J_{\alpha}\left({\sqrt{\beta t} \over t+p}\right)=
{\beta^{\alpha/2} \Gamma(\mu+1) \Gamma(\nu-\mu-1) \over 2^{\alpha}
p^{\nu-\mu-1} \Gamma(\alpha+1) \Gamma(\nu)}\;\;\times
&\nonumber\\&
_2F_3\left(\mu+1,\nu-\mu-1;{\nu \over 2},{\nu +1 \over 2},\alpha+1;
-{\beta \over 16 p}\right)\;\;,
\label{oneeightd}&\end{eqnarray}
\begin{eqnarray}&
\int_0^1x^{\mu}(1-x)^{\nu-\mu-2}e^{-\beta x^2}dx=B(\nu-\mu-1,\mu+1)
\;_2F_2(\frac{\mu+1}{2},\frac{\mu}{2}+1;\frac{\nu}{2},\frac{\nu+1}{2};
-\beta)
&\end{eqnarray}
\begin{eqnarray}&
\int_0^1x^{\nu-\mu-1}(1-x)^{\mu-1}e^{-zx^3}dx=\frac{\Gamma(\mu)
\Gamma(\nu-\mu)}{\Gamma(\nu)}
\;\;\times
&\nonumber\\&
_3F_3\left({\nu-\mu \over 3},{\nu-\mu+1 \over 3},{\nu-\mu+2 \over 3};
{\nu \over 3},{\nu+1 \over 3},{\nu +2 \over 3};-z\right),
\label{oneeightf}&\end{eqnarray}
\begin{eqnarray}&
\int_0^1x^{\frac{3}{2}\nu+1}(1-x)^{\mu-1}J_{\nu}(2zx^{3/2})dx=
&\nonumber\\&
\frac{\Gamma(\mu)\Gamma(3\nu+2)z^{\nu}}{\Gamma(3\nu-\mu+4)\Gamma(
\nu+1)}\;_2F_3(\nu+\frac{2}{3},\nu+\frac{4}{3};\nu+\frac{\mu+2}{3},
\nu+\frac{\mu+3}{3},\nu+\frac{\mu+4}{3};-z^2)
\label{oneeightg}&\end{eqnarray}
\begin{eqnarray}&
\int_0^1x^{\nu-1}(1-x)^{\mu-1}\;_1F_2(\frac{\nu+\mu}{2};\nu,
\frac{\nu-\mu-1}{2};z^2x^2)=
&\nonumber\\&
\frac{\sqrt{\pi}}{z^{\nu-1}2^{\mu}}\frac{\Gamma(\nu)\Gamma(\mu)
\Gamma(\frac{\nu-\mu+1}{2})}{\Gamma(\frac{\nu+\mu}{2})}
J_{\frac{\nu+\mu-1}{2}}(z)J_{\frac{\nu-\mu-1}{2}}(z).
\label{oneeighth}&\end{eqnarray}

There are many inverse LTs listed in Ref.\ \cite{pru92a}, which
have the form $x^{\alpha} (x+a)^{-\nu}$. For example, if we put
$\alpha= -\nu$ in Eq.\ (\ref{six}) and use our theorem, then for
${\rm Re}\;\nu>0$ we find that, expressed as a Mellin transform,
\begin{eqnarray}&
\int_0^{\infty}\frac{x^{s-1}}{[(x+a)(x+b)]^{\nu}}dx=
&\nonumber\\&
(\frac{2}{a+b})^{2\nu-s}B(s,2\nu-s)\;_2F_1(\nu-\frac{s}{2},
\nu-\frac{s-1}{2};\nu+\frac{1}{2};(\frac{a-b}{a+b})^2),
&\nonumber\\&
\label{onenine}&\end{eqnarray}
where we have used No. 2.15.3.2 from Ref.\ \cite{pru86}.
One can come up with many interesting
special cases, such as 
\begin{equation}
\int_0^{\infty}[\frac{x^2}{(x^2+a^2)(x^2+b^2)}]^{\nu}dx=
\frac{4^{\nu-1}\Gamma(\nu+1/2)\Gamma(\nu-1/2)}{\Gamma(2\nu)(a+b)^{
2\nu-1}}.
\label{twozero}\end{equation}

To  consider an example that does not appear in the standard
tables of integrals \cite{gra80}, \cite{ape83}, \cite{pru86a},
let
\begin{equation}
I_{\mu}(p)= \int_p^{\infty} dx\; {x^{-\alpha}\;(x^2+a^2)^{-\nu}
\over (x-p)^{\mu}} \quad.
\label{twotwo}\end{equation}
Then utilising our theorem we find that
\begin{equation}
I_{\mu}(p)= {\Gamma(1-\mu) \over \Gamma(\alpha+2\nu)}\;
\int_0^{\infty}dt\; t^{\mu+\alpha+2\nu-2} e^{-pt}\;
_1F_2\left(\nu;\nu+{\alpha \over 2},\nu+{\alpha+1 \over 2};-{a^2 t^2 \over 4}
\right)\quad ,
\label{twothree}\end{equation}
which is valid for Re $\alpha+2\nu >0$ and for Im $|a|=0$
according to p. 28 of Ref.\ \cite{pru92a}. Utilising No. 
3.38.1.17 of Ref.\ \cite{pru92}, we find that the above 
integral yields
\begin{eqnarray}&
I_{\mu}(p)={B(1-\mu,\mu+\alpha+2\nu-1) \over
p^{\mu+\alpha+2\nu-1}}\;\;\times
&\nonumber\\&
_3F_2\left(\nu+ {\alpha+\mu-1 \over 2},\nu+{\alpha+\mu \over 2},\nu;
\nu+{\alpha \over 2},\nu+{\alpha+1 \over 2};-{a^2 \over p^2}\right) \quad,
\label{twofour}&\end{eqnarray}
where Re $(\mu+\alpha+2\nu) > 1$.

For the case of $\alpha=0$, the above results yield the Mellin
transform
\begin{eqnarray}&
\int_0^{\infty}\frac{x^{s-1}}{[(x+1)^2+a^2]^{\nu}}dx=
&\nonumber\\&
B(s,2\nu-s)\;_2F_1(\nu-\frac{s}{2},\nu+\frac{1-s}{2};\nu+\frac{1}{2};-a^2),
\label{twofive}&\end{eqnarray}
which is valid for Re $0<Re\;s<2Re\;\nu$, and
\begin{eqnarray}&
\int_0^{\infty}\frac{x^{\nu-1/2}}{(x+1)[(x+1)^2+a^2]^{\nu}}dx=
&\nonumber\\&
(\frac{2}{a^2})^{\nu}\frac{\Gamma^2(\nu+1/2)}{\Gamma(\nu+1)}
(\sqrt{1+a^2}-1)^{\nu}\;_2F_1(\nu,\frac{1}{2};\nu+1;\frac{1+
\sqrt{1+a^2}}{2}).
\label{twoeight}&\end{eqnarray}
Many other interesting results can be obtained in this way; several are
given in appendix A.

Now if we return to Eq.\ (\ref{twofive}) and put $\nu=k+\nu$,
then after using the duplication formula for the gamma function,
we get
\begin{eqnarray}&
{\Gamma(k+\nu) \over \Gamma(k+\nu+\mu/2)}\int_p^{\infty} dx\;
{(x^2+a^2)^{-(k+\nu)} \over (x-p)^{\mu}}={2^{\mu-1} \Gamma(1-\mu)
\over (a^2+p^2)^{k+\nu+\mu/2-1/2}}\;{\Gamma(k+\nu+\mu/2-1/2) \over
\Gamma(k+\nu+1/2)}\times
&\nonumber\\&
_2F_1\left(k+\nu+{\mu-1 \over 2},{1-\mu \over 2};k+\nu+{1 \over 2};
{a^2 \over (a^2+p^2)}\right) \quad.
\label{threetwo}&\end{eqnarray}
Multiplying both sides of the above equation by $(\gamma)_k a^{2k}/k!$,
where $(\gamma)_k=\Gamma(k+\gamma)/\Gamma(\gamma)$ and then
summing from $k=0$ to $\infty$, we eventually arrive at
\begin{eqnarray}&
\int_p^{\infty}dx\; {(x^2+a^2)^{-\nu} \over (x-p)^{\mu}}\;
_2F_1\left(\gamma,\nu;\nu+{\mu \over 2};{a^2 \over x^2+a^2}\right)=
{\Gamma(2\nu+\mu-1)\;\Gamma(1-\mu) \over (a^2+p^2)^{\nu+\mu/2-1/2}\;
\Gamma(2\nu)}\;\;\times
&\nonumber\\&
_2F_1\left(\nu+{\mu-1 \over 2},\gamma+{1-\mu \over 2};\nu+{1 \over 2};
{a^2 \over a^2+p^2}\right)\;\;. 
\label{threethree}&\end{eqnarray}
In obtaining Eq.\ (\ref{threethree}) we have used No. 6.7.1.7
from Prudnikov et al \cite{pru90}. We may also multiply Eq.\ 
(\ref{twofive}) with $\nu=k+\nu$ by $t^k/k!\Gamma(k+\nu+1/2)$ 
and sum to obtain the Mellin transform
\begin{eqnarray}&
\int_0^{\infty}\frac{x^{s-1}}{[(x+1)^2+a^2]^{\nu}}\;_2F_1(
\frac{\nu}{2},\frac{\nu-1}{2};\nu+\frac{1}{2};\frac{b^2}{(x+1)^2})=
&\nonumber\\&
2^{-s}\Gamma(s)\;_2F_1(\nu-\frac{s}{2},\nu+\frac{1-s}{2};\nu+\frac{1}{2};
b^2-a^2)
\label{threefour}&\end{eqnarray}
where $|a| <1$, $|b| <1$, $|b^2-a^2| < 1$ and we have
used No. 6.7.1.11 from Ref.\ \cite{pru90}.
 
Before considering a numerical example, let us now consider the
following integral
\begin{equation}
I= \int_{0}^{1} dy\;{y^{\alpha-1} \over (1-y)^{\mu}}\;
_2F_1(a,b;c;-\omega y^2/p^2) \quad ,
\label{threefive}\end{equation}
where ${\rm Im}\; \omega=0$. By making the change of variable, $y=p/x$,
we find that Eq.\ (\ref{threefive}) becomes
\begin{equation}
I=p^{\alpha} \int_p^{\infty}dx\;{x^{\mu-\alpha-1} \over (x-p)^{\mu}}
\;_2F_1(a,b;c;-\omega/x^2) \quad,
\label{threesix}\end{equation}
whereupon we notice that the integral is now in the form of a fractional
integral. Thus after applying our theorem we find for
${\rm Re}\;\alpha>0$ that
\begin{eqnarray}&
I={p^{\alpha} \over \Gamma(\alpha+1-\mu)}\;\int_{0}^{\infty}dt\;
t^{\alpha-1}\;e^{-pt}\; _2F_3\left(a,b;c,{\alpha+1-\mu \over 2},
{\alpha+2-\mu \over 2};{-\omega t^2 \over 4} \right)\;\;=
&\nonumber\\&
{\Gamma(\alpha) \over \Gamma(\alpha+1-\mu)}\;
_4F_3\left(a,b,{\alpha \over 2},{\alpha +1 \over 2};c,
{\alpha +1-\mu \over 2},{\alpha +2-\mu \over 2};-{\omega \over p^2}\right)
\quad .
\label{threeseven}&\end{eqnarray}

We can use the above result in conjunction with our theorem and
specific inverse LTs to develop some interesting results. For
example, if we put $b=a+1$, $c=2a$ and $\alpha=\mu+2a$, then
Eq.\ (\ref{threefive}) can be written with the aid of our theorem
as
\begin{eqnarray}&
\int_0^1 dy\;{y^{\mu+2a-1} \over (1-y)^{\mu}}\;
_2F_1\left(a,a+1;2a;-{\omega y^2 \over p^2}\right)= {p^{\mu+2a} \Gamma^2(a+1/2)
\over 2^{2-4a}\omega^{a-1/2}\Gamma(2a+1)}\int_0^{\infty}dt\;t^{\mu}\;\times
&\nonumber\\&
e^{-pt} J_{a-1/2}^2(\sqrt{\omega}t/2)\;=
{\Gamma(\mu+2a) \over \Gamma(2a+1)}\;
_3F_2\left(a,a+{\mu \over 2},a+{\mu+1 \over 2};2a,a+{1 \over 2};-{\omega
\over p^2}\right) \;\;,
\label{threeeight}&\end{eqnarray}
where ${\rm Re}\; a >1/2$ and we have used No. 3.35.1.20 from Ref.
\cite{pru92a}. The above is basically the LT of $t^{\mu}
J_{a-1/2}^2(\omega^{1/2}t/2)$. Hence, our theorem can be used to
develop new LTs as well as Mellin transforms.

To complete this section we now consider the following numerical example
\begin{equation}
I_{\mu}(p)= \int_0^{\infty}dx\; {\ln(x^2+a^2) \over (x^2+a^2)^{1/2}
(x-p)^{\mu}} \quad ,
\label{fourthree}\end{equation}
where Im $a=0$. Then by using our theorem we get
\begin{equation}
I_{\mu}(p)= -\Gamma(1-\mu)\int_0^{\infty} dt\; t^{\mu-1}e^{-pt}\;\left(
{\pi \over 2}Y_0(at)+\left(C+\ln{2t \over a}\right)J_0(at)\right)\quad,
\label{fourfour}\end{equation}
where $C$ denotes Euler's constant, $0 <$Re $\mu <1$ and $Y_0(at)$ is the
Neumann function. We can evaluate the above integral, but the final form
is not very useful as we shall see. By using Nos. 2.12.8.4 and 2.13.6.3
from Ref.\ \cite{pru86} the above integral yields after a little
manipulation
\begin{eqnarray}&
I_{\mu}(p)={\Gamma(\mu) \Gamma(1-\mu) \over (p^2+a^2)^{\mu/2}}\;
\biggl[(\ln (a/2)-C)
P_{\mu-1}\left({p \over (p^2+a^2)^{1/2}}\right)+Q_{\mu-1}\left(
{p \over (p^2+a^2)^{1/2}}\right)\biggr]\;\;-
&\nonumber\\&
\Gamma(1-\mu){\partial \over \partial \alpha}\left({\Gamma(\mu+\alpha) \over
(p^2+a^2)^{(\mu+\alpha)/2}} P_{\mu+\alpha-1}\left({p \over
(p^2+a^2)^{1/2}}\right)\right)_{\alpha=0} \quad,
\label{fourfive}&\end{eqnarray}
where $P_{\nu}(x)$ and $Q_{\nu}(x)$ are respectively the Legendre functions
of the first and second kind. As can be seen the above result is
awkward because of the partial derivative. Since $J_0(at)$ is bounded
and $Y_0(at)$ can be written as the sum of a bounded term and
$2\ln(at/2)$, Eq.\ (\ref{fourfour}) is simply the sum of an integral
yielding the gamma function and another integral of the form
\begin{equation}
I=\int_0^{\infty}dt\;t^{\mu-1}\;e^{-pt} \ln\;t \quad,
\label{foursix}\end{equation}
which is rapidly converging at the upper end. The integral can
be handled at the lower end provided Re $\mu >0$ because then the
algebraico-logarithmic singularity can be numerically integrated
by following the procedure as described by Piessens et al on
p. 47 of Ref.\ \cite{pie83}.

\section{Non-fractional Integrals}
In the previous section we remarked that our theorem could be
extended to integrals not necessarily of a fractional form. 
To see this,
note that we can write any integral as 
\begin{equation}
I=\int_0^{\infty}dx\;F_1(x) F_2(x) \quad.
\label{fourseven}\end{equation}
We should add at this stage that any integral with other definite
limits can be converted to one with the above limits. Assuming that
the inverse LT of $F_1(x)$ exists, we can write Eq.\ (\ref{fourseven})
as
\begin{eqnarray}&
I=\int_0^{\infty}dx\;F_2(x) \int_0^{\infty}dt\;e^{-xt}\;
{\cal L}^{-1}_t(F_1)= \int_0^{\infty}dt\;{\cal L}_t(F_2)
{\cal L}^{-1}_t(F_1) \quad,
\label{foureight}&\end{eqnarray}
where we have assumed that it is valid to change the order of
integration. Now one can see that if the inverse LT of $F_1(x)$
has an exponential factor in it as, for example, the inverse LT of
$(x-p)^{-\mu} \Theta(x-p)$ has, then Eq.\ (\ref{foureight}) simply
becomes another LT. However, this is not necessary
, although such a factor is particularly advantageous
in view of (1) the availability of tables of LTs and (2) its
amenability to numerical computation. As a consequence, in this section
we shall concentrate on inverse LTs of $F_1(x)$ that have some form of
decaying exponential behaviour, e.g. the MacDonald function, $K_{\nu}(x)$.

Let us begin by considering a simple example, which does not appear
in the standard tables of integrals \cite{gra80}, \cite{ape83},
\cite{pru86a} but is a special case of an example considered in the
previous section:
\begin{equation}
I=\int_p^{\infty}dx\;{x^{\nu+1}(x-p)^{\nu} \over (x^2+a^2)^{3/2}}
\quad,
\label{fournine}\end{equation}
where the integral is defined for $-1<$ Re $\nu<1/2$. Since
$ {\cal L}^{-1}(x/(x^2+a^2)^{3/2})=t J_0(at)$, the above
integral can be written as
\begin{eqnarray}&
I=\int_0^{\infty}dt\;t \;J_0(at) {\cal L}_t(x^{\nu}(x-p)^{\nu})\;=
{\Gamma(\nu+1)\; p^{\nu+1/2} \over \pi^{1/2}}\;\;\times
&\nonumber\\&
\int_0^{\infty}dt\;t^{-\nu+1/2}\;e^{-pt/2}\;J_0(at)K_{\nu+1/2}(pt/2)\quad.
\label{fivezero}&\end{eqnarray}
Eq.\ (\ref{fivezero}) is still unknown but can be evaluated by
expanding the Bessel function into a series. Hence, by
utilising No. 6.621.3 in Ref.\ \cite{gra80}, we get
\begin{equation}
I={\Gamma(\nu+1)\;\Gamma(1-2\nu) \over p^{1-2\nu}\;\Gamma(2-\nu)}
\;_3F_2\left({1 \over 2}-\nu,1-\nu,{3 \over 2};1-{\nu \over 2},
{3-\nu \over 2};-{a^2 \over p^2}\right)\quad.
\label{fiveone}\end{equation}
Eq.\ (\ref{fiveone}) reduces to the known result of $(a^2+p^2)^{-1/2}$
for $\nu=0$. The result can also be checked by putting 
$\mu=-\nu$, $\nu=3/2$ and $\alpha=-\nu-1$ into Eq.\ (\ref{twofour}).

We can use the modification of our theorem to display
similarities in apparently different integrals. For example,
consider the following integrals:
\begin{eqnarray}&
I_1=\int_a^{\infty}dx\;x^{-\mu}e^{-\beta/x}\;\Bigl[(x+\sqrt{x^2-a^2})^{\nu}
+ (x-\sqrt{x^2-a^2})^{\nu}\Bigr](x^2-a^2)^{-1/2} \quad,
\label{fivefive}&\end{eqnarray}
and
\begin{eqnarray}&
I_2=\int_a^{\infty}dx\;x^{-\mu}e^{-\beta/x}\;(x^2-a^2)^{\nu}
\quad.
\label{fivesix}&\end{eqnarray}
For $I_1$ we require Re $\mu >$ Re $\nu +1$ while for $I_2$, Re $\mu >
2$ Re $\nu +1$ and Re $\nu > -1$. To evaluate both integrals we also require
${\cal L}_t^{-1}(x^{-\mu}e^{-\beta /x})=(t/\beta)^{(\mu-1)/2}
J_{\mu-1}(2\sqrt{\beta t})$. Then  we have
\begin{equation}
I_1=2a^{\nu}\beta^{(1-\mu)/2} \int_0^{\infty}dt\;t^{(\mu-1)/2}
J_{\mu-1}(2\sqrt{\beta t}) K_{\nu+1/2}(at) \quad,
\label{fiveseven}\end{equation}
and
\begin{equation}
I_2={(2a)^{\nu+1/2}\; \Gamma(\nu+1) \over  \beta^{(\mu-1)/2}\;
\sqrt{\pi}}\int_0^{\infty}dt\;t^{\mu/2-\nu-1}
J_{\mu-1}(2\sqrt{\beta t}) K_{\nu}(at) \quad.
\label{fiveeight}\end{equation}
Now we can see that both integrals are of the same type, which is
given as No. 2.16.22.6 in Ref.\ \cite{pru86}, and find that
\begin{eqnarray}&
I_1=2^{\mu-1} a^{\nu-\mu-1}\biggl[a\;B\left({\mu+\nu \over 2}+{1
\over 4},{\mu-\nu \over 2}-{1 \over 4}\right)\;\;\times
&\nonumber\\&
_2F_3\left({\mu+\nu \over 2}+{1 \over 4},{\mu-\nu \over 2}-{1 \over 4};
{1 \over 2},{\mu \over 2},{\mu +1 \over 2};{\beta^2 \over 4a^2}\right) -
2 \beta \;B\left({\mu+\nu \over 2}+{3 \over 4},{\mu-\nu \over 2}+
{1 \over 4}\right)\;\times
&\nonumber\\&
 _2F_3\left({\mu+\nu \over 2}+{3 \over 4},{\mu-\nu \over
2}+{1 \over 4};{3 \over 2}, {\mu+1 \over 2},{\mu +2\over 2};
{\beta^2 \over 4a^2}\right)\biggr]\;,
\label{fivenine}&\end{eqnarray}
which is valid for Re $a>0$ and ${\rm Re}\; \mu > |{\rm Re}(\nu+1/2)|$.
Eq.\ (\ref{fiveeight}) becomes
\begin{eqnarray}&
I_2={2^{\mu-2} a^{2\nu-\mu} \Gamma(\nu+1) \over \sqrt{\pi}}\;
\biggl[{a\; \Gamma(\mu/2-1/4)\; \Gamma(\mu/2-\nu-1/4) \over \Gamma(\mu)}
\;\;\times
&\nonumber\\&
_2F_3\left({\mu \over 2}-{1 \over 4},{\mu \over 2}-\nu-{1 \over 4};
{1 \over 2},{\mu \over 2},{\mu +1 \over 2};{\beta^2 \over 4a^2}\right)
\;- {2 \beta \;\Gamma(\mu/2+1/4)\; \Gamma(\mu/2-\nu+1/4) \over \Gamma(\mu+1)}
\;\;\times
&\nonumber\\&
_2F_3\left({\mu \over 2}+{1 \over 4},{\mu \over 2}-\nu+ {1 \over 4};
{3 \over 2},{\mu +1 \over 2}, {\mu \over 2}+1;{\beta^2 \over 4a^2}\right)
\biggr] \quad,
\label{sixzero}&\end{eqnarray}
the latter integral being valid for Re $\mu >|$ Re $ \nu| + $ Re $ \nu
+1/2$ and Re $a>0$.

We should also point out that occasionally a fractional
integral, which may not solved by using our theorem can be solved
via the modified approach presented in this section. As an example, let
us consider the following fractional integral
\begin{equation}
I=\int_p^{\infty} dx\;(x-p)^{\mu}x^{\mu}\;(x^2-px+b)^{-1/2}
\quad,
\label{sixone}\end{equation}
which is defined for $-1 < $ Re $\mu < 0$. To evaluate the above
integral we shall also require that $|$ Im $(b-p^2/4)^{1/2}| - $ Re $p/2 < 0$.
In order to utilise the theorem given in the previous section we
require the inverse LT of $x^{\mu}(x^2-px+b)^{-1/2}$, which is
not given in Refs. \cite{obe73} and \cite{pru92a}. However, the inverse
LT of $(x^2-px+b)^{-1/2}$ appears as No. 2.1.6.1 in Prudnikov et al
\cite{pru92a} and hence, modification of our theorem yields 
\begin{eqnarray}&
I={p^{\mu+1/2} \;\Gamma(\mu+1) \over \sqrt{\pi}}\int_0^{\infty}dt\;
t^{-(\mu+1/2)}\;J_0(t\sqrt{b-p^2/4})\;K_{\mu+1/2}(pt/2)\;\;=
&\nonumber\\&
{p^{2\mu} \over 2^{2\mu+1}}\;\Gamma(\mu+1)\;\Gamma(-\mu)\;
_2F_1\left({1 \over 2},-\mu;1;1-{4b \over p^2}\right) \quad. 
\label{sixtwo}&\end{eqnarray}
In obtaining Eq.\ (\ref{sixtwo}) we have used No. 2.16.21.1 from Ref.\
\cite{pru86}. For the special case of $\mu=-1/2$, we find that 
$I=2p^{-1}{\bf K}((1-4b/p^2)^{1/2})$ where ${\bf K}(x)$ is the complete
elliptic integral of the first kind. Finally, more results based on the
modification of our theorem appear in Appendix B.

\section{Further Extensions}
In this section we consider some extensions of our study of
fractional integrals, the first dealing with integrals of the 
form of
\begin{equation}
I=\int_{p}^{\infty}dx\;(x-p)^{-\mu}\;(x-p_1)^{-\mu_1}\;F(x) \;,
\label{sixthree}\end{equation}
where $p >p_1$, Re $\mu_1 >0$ and $0 < $ Re $\mu <1$. Utilising the
integral representation for the gamma function, we can write Eq.\
(\ref{sixthree}) as
\begin{equation}
I={1 \over \Gamma(\mu) \Gamma(\mu_1)} \int_{p}^{\infty}dx\;
\int_0^{\infty}dt\;e^{-xt}\;G(t)\;\int_0^{\infty}dz\;
z^{\mu-1}\;e^{-(x-p)z}\int_0^{\infty}dy\;y^{\mu_1-1}\;e^{-(x-p_1)y}
\quad ,
\label{sixfour}\end{equation}
where $G(t)$ is the inverse LT of $F(x)$. After a change of variable
and evaluation of some of the integrals, one arrives at
\begin{eqnarray}&
I={\Gamma(1-\mu) \over \Gamma(\mu_1)}\;\int_0^{\infty}dt\;e^{-pt}\;
G(t)\;\int_0^{\infty}dy\;y^{\mu_1-1}\;e^{-(p-p_1)y}\;(y+t)^{\mu-1}\;=
\Gamma(1-\mu)\;\;\times
&\nonumber\\&
(p-p_1)^{-(\mu+\mu_1)/2}\;\int_0^{\infty}dt\;
e^{-(p+p_1)t/2}\;G(t) t^{(\mu+\mu_1)/2-1}\;W_{(\mu-\mu_1)/2,
(\mu+\mu_1-1)/2}\Bigl((p-p_1)t\Bigr)\;,
\label{sixfive}&\end{eqnarray}
where we have used No.\ 3.383.8 from Gradshteyn and Ryzhik \cite{gra80}
and $W_{\nu,\mu}(z)$ denotes the Whittaker function. This result can
be analytically continued to ${\rm Re}\;\mu<1$ for any value of $\mu_1$
provided the integrals in Eqs. (\ref{sixthree}) and (\ref{sixfive})
are defined. The condition on ${\rm Re}\;\mu$ arises from the factor
of $\Gamma(1-\mu)$ in Eq.\ (\ref{sixfive}). If we had introduced the
integral representation for the gamma function in our theorem in Sec.\
2, then we would have found that $0<{\rm Re}\;\mu<1$ to obtain Eq.\
(\ref{three}). However, we have shown that our theorem is valid for
${\rm Re}\;\mu<1$ by using LT theory. This condition on $\mu$ arises
due to the $\Gamma(1-\mu)$ factor in Eq.\ (\ref{three}) but there
is no lower bound on $\mu$. Thus, while using the integral
representation for the gamma function is initially restrictive, analytic
continuation can be used to extend results such as Eqs. (\ref{three}) and
(\ref{sixfive}) beyond these restrictions.

As an example of Eq.\ (\ref{sixfive}) in action, let us put
$F(x)=(x^2+ax+b)^{-1/2}$, $\mu_1=1$ and $\mu=0$. Then the original
integral given by Eq.\ (\ref{sixthree}) can be evaluated by using No.
2.2.9.41 of Ref.\ \cite{pru86a}. The inverse LT of $F(x)$ is given as
No. 2.1.6.1 in Ref.\ \cite{pru92a}. Hence from Eq.\ (\ref{sixfive}),
we obtain
\begin{eqnarray}&
I=\int_0^{\infty}dt\;t^{-1/2}\;e^{-(p+p_1+a)t/2} J_0(t\sqrt{b^2-a^2/4})
W_{-1/2,0}((p-p_1)t)\;=
&\nonumber\\&
\sqrt{{p-p_1 \over p_1^2+ap_1+b}}\;\ln\;{(p^2+ap+b)^{1/2}
(p_1^2+ap_1+b)^{1/2}+b+ap/2+p_1(p+a/2) \over (p-p_1)(p_1+a/2
+(p_1^2+ap_1+b)^{1/2})}\;\;,
\label{sixfivea}&\end{eqnarray}
where $(p_1^2+ap_1+b)^{1/2}>0$, $(p^2+ap+b)^{1/2}>0$,
$(p+a/2)>-(p^2+ap+b)^{1/2}$ and ${\rm Re}\;p>
|{\rm Im}(b^2-a^2/4)^{1/2}|-{\rm Re}\;a/2$. For $(-p_1^2-ap_1-b)^{1/2}>0$ 
and $(p+a/2) > (p^2+ap+b)^{1/2}$, we find that 
\begin{eqnarray}&
I=\sqrt{{p_1-p \over (p_1^2+ap_1+b)}}\;\arccos\left(1+
{p_1^2+ap_1+b \over (p-p_1)(\sqrt{p^2+ap+b}+p+a/2)}\right)\;\;.
\label{sixfiveb}&\end{eqnarray}

We can also develop another integral representation for Eq.\ (\ref{sixthree})
by using the Feynman integral \cite{sch61}, \cite{pok89} which allows
us to replace the denominator of Eq.\ (\ref{sixthree}) by
\begin{equation}
(x-p)^{-\mu}\;(x-p_1)^{-\mu_1}={\Gamma(\mu+\mu_1) \over \Gamma(\mu)
\Gamma(\mu_1)}\;\int_0^{1}dt\;{t^{\mu_1-1}\;(1-t)^{\mu-1} \over
((p-p_1)t+x-p)^{\mu+\mu_1}} \quad,
\label{sixsix}\end{equation}
where ${\rm Re}\;\mu,\mu_1>0$. In Appendix C we establish the above result
from its more general version, which appears as No. 2.2.6.1 in Ref.\
\cite{pru86a} and then use these integrals to develop several interesting
results. Thus, Eq.\ (\ref{sixthree}) can be written as
\begin{eqnarray}&
I={\Gamma(\mu+\mu_1) \over \Gamma(\mu) \Gamma(\mu_1)} \int_0^{\infty}
dt\;e^{-pt}G(t) \int_0^1ds\;s^{\mu_1-1}(1-s)^{\mu-1} \int_0^{\infty}
dy\;{e^{-yt} \over (y+(p-p_1)s)^{\mu+\mu_1}}\quad.
\label{sixseven}&\end{eqnarray}
The final integral in Eq.\ (\ref{sixseven}) can also be evaluated by using
No. 3.383.4 in Gradshteyn and Ryzhik \cite{gra80}. So, Eq.\ (\ref{sixseven})
becomes
\begin{eqnarray}&
I={\Gamma(\mu+\mu_1) \over \Gamma(\mu) \Gamma(\mu_1)}\int_0^{\infty}dt\;
e^{-pt} G(t) t^{(\mu+\mu_1-2)/2}\int_0^1ds\;s^{\mu_1-1}(1-s)^{\mu-1}
((p-p_1)s)^{-(\mu+\mu_1)/2}\;\;\times
&\nonumber\\&
e^{(p-p_1)st/2} W_{-(\mu+\mu_1)/2,(1-\mu-\mu_1)/2}
\Bigl((p-p_1)st\Bigr)\quad .
\label{sixeight}&\end{eqnarray}
If in Eqs. ({\ref{sixfive}) and (\ref{sixeight}) we put $G(t)=J_0(zt)$,
i.e. $F(x)=(x^2+z^2)^{-1/2}$, then by multiplying both equations by
$zJ_0(zy)$ and integrating over $z$ from zero to infinity, we obtain
the interesting result of
\begin{eqnarray}&
\int_0^1dt\;t^{(\mu_1-\mu)/2-1}(1-t)^{\mu-1}e^{(p_1-p)t/2}\;
W_{(\mu_1+\mu)/2,(1-\mu-\mu_1)/2}\Bigl((p-p_1)t\Bigr)={\Gamma(1-\mu)
\Gamma(\mu) \over \Gamma(\mu+\mu_1)}\;\;\times
&\nonumber\\&
\Gamma(\mu_1)\;e^{(p-p_1)/2}\;W_{(\mu-\mu_1)/2,(\mu_1+\mu-1)/2}
\Bigl(p-p_1\Bigr)\quad,
\label{sixnine}&\end{eqnarray}
where $0< {\rm Re}\;\mu<1$ and ${\rm Re}\;\mu_1>0$.

For the special case of $\mu=\mu_1$ the integral given by Eq.\
(\ref{sixthree}) reduces to
\begin{eqnarray}&
\int_p^{\infty}dx\;(x^2-(p+p_1)x+pp_1)^{-\mu} F(x)=
{\Gamma(1-\mu) \over \sqrt{\pi}\;(p-p_1)^{\mu-1/2}}\;
\int_0^{\infty}dt\;e^{-(p+p_1)t/2}G(t)\;\;\times
&\nonumber\\&
t^{\mu-1/2}\;K_{\mu-1/2}\left((p-p_1)t/2\right)\quad,
\label{sevenzero}&\end{eqnarray}
where we have used Eq.\ (\ref{sixfive}) and the fact that
$W_{0,\mu}(z)=\sqrt{z/\pi}K_{\mu}(z/2)$ according to No.\ 9.235.2
of Ref.\ \cite{gra80}. Eq.\ (\ref{sevenzero}) can be checked by putting
$F(x)= x^{-\nu}$ whereupon the integral on the lhs after a change of
variable becomes
\begin{eqnarray}&
I=(p-p_1)^{-\mu} p^{1-\mu-\nu} \int_0^{\infty}ds\;s^{-\mu}
\left(1+{ps \over p-p_1}\right)^{-\mu}(1+s)^{-\nu}\;=
&\nonumber\\&
(p-p_1)^{-\mu}p^{1-\mu-\nu} B(1-\mu,\nu+2\mu-1)\;_2F_1\left(\mu,
1-\mu;\mu+\nu;-{p_1 \over p-p_1}\right)\quad.
\label{sevenzeroa}&\end{eqnarray}
In obtaining the above integral we have used No.\ 3.197.5 from Ref.
\cite{gra80}. It is valid for ${\rm Re}\;\mu<1$, $|{\rm arg}(p/(p-p_1))|<\pi$
and ${\rm Re}\;\nu>1-2{\rm Re}\;\mu$. One can show that the rhs
of Eq.\ (\ref{sevenzero}) gives the same result by noting that 
${\cal L}^{-1}_{t}(x)$ is $t^{\nu-1}/\Gamma(\nu)$ and utilising
Nos. 6.621.3 and 9.131.1 from Ref.\ \cite{gra80}. Furthermore, if we
put $p_1=-p$ in Eq.\ ({\ref{sevenzero}), then we obtain
\begin{eqnarray}&
\int_p^{\infty}dx\;(x^2-p^2)^{-\mu} F(x)={\Gamma(1-\mu) \over 
\sqrt{\pi}}\int_0^{\infty}dt\;(t/2p)^{\mu-1/2}\;G(t) K_{\mu-1/2}(pt)
\quad,
\label{sevenone}&\end{eqnarray}
which can also be found by following the procedure outlined in the
previous section where one would eventually evaluate the LT of
$(x^2-p^2)^{\mu}$. Using the latter approach, however, means that the
restriction on $\mu$ as given under Eq.\ (\ref{sixthree}) can be
analytically continued to ${\rm Re}\;\mu<1$.

A pleasant example that can be derived from Eq.\ (\ref{sevenzero}) is 
\begin{eqnarray}&
\int_0^{\infty}\frac{(x+c)}{(x+a)(x+b)}[x(x+1)]^{\nu}dx=
&\nonumber\\&
\frac{B(1+\nu,-2\nu)}{a-b}\{\frac{a-c}{2a-1}\;_2F_1(\frac{1}{2},
1;1-\nu;4\frac{a(a-1)+2c(1-c)}{(2a-1)^2})-
&\nonumber\\&
\frac{b-c}{2b-1}\;_2F_1(\frac{1}{2},1;1-\nu;4\frac{b(b-1)+2c(1-c)}{
(2b-1)^2}
\label{sevenoneb}&\end{eqnarray}
where, to ensure convergence, $-1<Re\;\nu<1$. In obtaining this we have utilised various results given
in Refs.\ \cite{pru86} and \cite{pru92a}.

To set the scene for the material to follow let us examine
Eq.\ (\ref{sevenone}) more closely by analysing the behaviour at the
lower and upper ends of the integrals on both sides. Clearly for the
integral on the lhs to be defined, $F(x)$ must be $O((x-p)^{\beta})$ as
$x \to p^{+}$ where ${\rm Re}\;\beta> {\rm Re}\;\mu-1$. This condition has
been implicitly assumed throughout in discussing all fractional integrals.
At the upper end, $F(x)$ must be $O(x^{\alpha})$, where 
${\rm Re}\;\alpha<2{\rm Re}\;\mu-1$. If we were to replace $\mu$ by
$\mu-n$, where $n$ is a positive integer or equal to zero, then the lhs
of Eq.\ ({\ref{sevenone}) would only be defined if
${\rm Re}\;\alpha<2{\rm Re}\;\mu-2n-1$. Thus, if we wanted to sum over
all values of $n$ from zero to infinity, then to guarantee convergence 
$F(x)$ would have to possess some exponential decaying factor such as
$\exp(-\epsilon x)$ as $x \to \infty$, where ${\rm Re}\; \epsilon>0$. However,
this may be too restrictive for many integrals may be defined after the
summation process has been effected.

The singular behaviour is reciprocated when examining the rhs of Eq.\
(\ref{sevenone}). Thus, if we were to replace $\mu$ by $\mu-n$ in the
rhs of this equation, then for the integral to be defined for $\mu-n-1/2<0$
we see that $G(t)$ must be $O(t^{\alpha})$ as $t \to 0^+$, where 
${\rm Re}\;\alpha>2n-2{\rm Re}\;\mu$. Again, we can overcome this singular
behaviour by requiring $G(t)$ to have an exponential factor. In this case
it would be $\exp(-\gamma/t)$ but it would also be too restrictive on $G(t)$.

By taking note of the above, we can introduce exponential factors to
both sides of Eq.\ (\ref{sevenone}) so that it becomes
\begin{eqnarray}&
\lim_{\epsilon \to 0^{+}}\;\Gamma(1-\mu)^{-1} \int_p^{\infty}dx\;
(x^2-p^2)^{-\mu} e^{-\epsilon (x-p)} F(x)=\lim_{\epsilon \to 0^{+}}\;
\int_0^{\infty}dt\;e^{-\epsilon/t} (t/2p)^{\mu-1/2}\times
&\nonumber\\&
\pi^{-1/2}\; G(t) K_{\mu-1/2}(pt)\quad.
\label{seventwo}&\end{eqnarray}
This enables us to replace $\mu$ by $\mu-n$ and then to sum from $n=0$
to $\infty$. As a consequence of absolute convergence, we can interchange the
order of the summations and integrations. Finally, we can take the limit as
$\epsilon \to 0^+$, which will yield integrals with less restrictions on
the behaviour of $F(x)$ and $G(t)$ than those mentioned previously. Thus
by putting $\mu=\mu-n$ and multiplying both sides by 
$(-1)^n \beta^{2n-\mu}/2^{2n-\mu}n!$ of Eq.\ (\ref{seventwo}), we obtain
\begin{eqnarray}&
\lim_{\epsilon \to 0^{+}}\;\sum_{n=0}^{\infty}{(-1)^n \beta^{2n-\mu} \over
2^{2n-\mu} n! \Gamma(1-\mu+n)}\;\int_p^{\infty}dx\;(x^2-p^2)^{-\mu+n}
e^{-\epsilon (x-p)}F(x)= \lim_{\epsilon \to 0^{+}} 
\sum_{n=0}^{\infty}{(-1)^n \over 2^{2n-\mu}}\;\times
&\nonumber\\& 
{\sqrt{\pi} \beta^{2n-\mu} \over  n!} \int_0^{\infty}dt\;e^{-\epsilon/t}
(t/2p)^{\mu-n-1/2}G(t) K_{n+1/2-\mu}(pt) \;\;.
\label{seventhree}&\end{eqnarray}
Since the integrals are now absolutely convergent, the order of the sums
and integrals can be interchanged to yield
\begin{eqnarray}&
\int_p^{\infty}dx\;(x^2-p^2)^{\mu/2} J_{\mu}(\beta \sqrt{x^2-p^2})
F(x)= \int_0^{\infty}dy\;{y^{\mu+1} \over
\sqrt{y^2+p^2}}\;J_{\mu}(\beta y) F(\sqrt{y^2+p^2})\;=
&\nonumber\\&
\sqrt{{2 \over \pi}}\;p^{\mu-1/2}\beta^{\mu}\int_0^{\infty}dt\;
(t^2+\beta^2)^{-\mu/2-1/4}G(t) K_{\mu+1/2}(p\sqrt{t^2+\beta^2})
\quad,
\label{sevenfour}&\end{eqnarray}
where we have taken the limit after evaluating the sums and have also
replaced $\mu$ by $-\mu$, so that the above result is valid for
${\rm Re}\;\mu>-1$. In addition, we have utilised the multiplication
theorem for the MacDonald function, which appears as No. 14 on
p. 112 of Ref.\ \cite{man65}.

Let us test the result given above by putting $F(x)=x^{-2\nu-1}$.
Then the lhs of Eq.\ ({\ref{sevenfour}) becomes No. 6.565.4 in Gradshteyn
and Ryzhik \cite{gra80}, which is
\begin{eqnarray}&
\int_0^{\infty}dy\;{y^{\mu+1} J_{\mu}(\beta y) \over (y^2+p^2)^{\nu+1}}
={p^{\mu-\nu} \beta^{\nu} \over 2^{\nu} \Gamma(\nu+1)}\;
K_{\mu-\nu}(\beta p) \quad.
\label{sevenfive}&\end{eqnarray}
This result is valid for $-1< {\rm Re}\;\mu<{\rm Re}\;(2\nu+3/2)$ with
both $\beta$ and $p$ greater than zero. According to No. 2.1.1.1 in
Ref.\ \cite{pru92a}, the inverse LT of $x^{-2\nu-1}$ is $t^{2\nu}/
\Gamma(2\nu+1)$, which when introduced into the rhs of Eq.\ (\ref{sevenfour})
yields the integral given by No. 6.596.3 of Gradshteyn and Ryzhik.
This, in turn, yields the rhs of Eq.\ (\ref{sevenfive}) after application
of the duplication formula for the gamma function and is valid for
${\rm Re}\;\nu>-1/2$. Thus the integral emanating from the rhs of 
Eq.\ (\ref{seventwo}) represents the analytic continuation of the integral
on the lhs of Eq.\ (\ref{sevenfour}).

We now consider an example, which does not appear in the standard
tables of integrals. This is
\begin{eqnarray}&
I=\int_0^{\infty}dt\;(t^2+\beta^2)^{-\mu/2-1/4}t^{4\mu}\;
_1F_2\left(\mu+{1 \over 2};2\mu+{1\over 2}, 2\mu+1;
{p^2t^2 \over 4}\right)K_{\mu+1/2}(p\sqrt{t^2+\beta^2}) \;\;.
\label{sevensix}&\end{eqnarray}
Utilising No. 2.1.7.20 of Ref.\ \cite{pru92a} and Eq.\ (\ref{sevenfour})
above, we can transform Eq.\ (\ref{sevensix}) into
\begin{eqnarray}&
I=\sqrt{{\pi \over 2}}\;{\Gamma(4\mu+1) \over \beta^{\mu} p^{\mu-1/2}}
\int_0^{\infty}dy\;y^{-\mu}(y^2+p^2)^{-\mu-1/2}\;
J_{\mu}(\beta y)=\sqrt{{\pi \over 2}}\;\Gamma(\mu+1)\;\;\times
&\nonumber\\&
{2^{\mu}\; \Gamma(4\mu+1) \over p^{3\mu-1/2}\; \Gamma(2\mu+1)}\;
I_{\mu}\left({\beta p \over 2}\right)K_{\mu}\left({\beta p \over
2}\right) \quad,
\label{sevenseven}&\end{eqnarray}
where ${\rm Re}\;\mu>-1/2$, $\beta>0$ and we have used No. 2.12.4.30 from
Ref. \cite{pru86}.

An interesting application can be found by putting $\mu=-1/2$ in Eq.\
(\ref{sevenfour}), which yields
\begin{eqnarray}&
\int_0^{\infty}dy\;{\cos(\beta y) \over \sqrt{y^2+p^2}}\;
F\left(\sqrt{y^2+p^2}\right)=p^{-1}\int_0^{\infty}dt\;G(t) 
K_0\left(p\sqrt{t^2+\beta^2}\right)\;\;.
\label{sevensevena}&\end{eqnarray}
If we put $G(t)=K_1(p(t^2+p^2)^{1/2})/(t^2+p^2)^{1/2}$, then using
the result in the appendix of Ref.\ \cite{bur86} we find that the
integral on the rhs of Eq.\ (\ref{sevensevena}) yields $\pi K_0(2p\beta)
/2\beta p^2$. We can identify $F((y^2+p^2)^{1/2})$ by noting from
No.\ 6.726.4 of Gradshteyn and Ryzhik \cite{gra80} that
\begin{eqnarray}&
\sqrt{{3\pi \over 2\beta}}\int_0^{\infty}dy\;{\cos(\beta y) \over
p^2(y^2+p^2)^{1/4}}\;K_{1/2}\left(\beta \sqrt{3(y^2+p^2)}\right)=
{\pi \over 2 \beta p^2}\;K_0(2p\beta)\;\;.
\label{sevensevenb}&\end{eqnarray}
Hence, we obtain the following LT which does not appear in Refs.
\cite{obe73} and \cite{pru92}:
\begin{eqnarray}&
\int_0^{\infty}dt\;e^{-yt}\;(t^2+\beta^2)^{1/2}K_{1}(p\sqrt{t^2+\beta^2})=
\sqrt{{3 \pi y \over 2 \beta}}\;p^{-2} K_{1/2}(\sqrt{3}\beta y)\;\;.
\label{sevensevenc}&\end{eqnarray}

We should mention that the approach used to derive Eq.\ (\ref{sevenfour})
from Eq.\ (\ref{sevenone}) can also be applied to Eq.\ (\ref{sevenzero}),
whereupon one obtains
\begin{eqnarray}&
\int_p^{\infty}dx\;(x-p)^{\mu/2}(x-p_1)^{\mu/2} J_{\mu}(\beta
\sqrt{(x-p)(x-p_1)}) F(x)={\beta^{\mu} (p-p_1)^{\mu+1/2} \over 2^{\mu}
\sqrt{\pi}}\;\;\times
&\nonumber\\&
\int_0^{\infty}dt\;e^{-(p+p_1)t/2}G(t)
(t^2+\beta^2)^{-\mu/2-1/4} K_{\mu+1/2}\left({(p-p_1) \over 2} 
\sqrt{t^2+\beta^2}\right)\quad.
\label{seveneight}&\end{eqnarray}
In the above result we have put $\mu=-\mu$, so that it is valid for 
${\rm Re}\;\mu>-1$. Eq.\ (\ref{seveneight}) can be verified by multiplying
both sides by $\beta J_{\mu}(\beta y)$ and then integrating over $\beta$ from
zero to $\infty$. This gives the correct result of $F(x)$ as the LT of
$G(t)$ after utilising No. 6.596.7 from Gradshteyn and Ryzhik \cite{gra80}.  
In addition, putting $p_1=-p$ yields Eq.\ (\ref{sevenfour}). By making
an appropriate change of variable, one can express Eq.\ (\ref{seveneight})
in the more convenient form of
\begin{eqnarray}&
\int_0^{\infty}ds\;{s^{\mu+1}\;J_{\mu}(\beta s) \over \sqrt{\alpha^2+
4s^2}}\;F\left({\sqrt{\alpha^2+4s^2}+\alpha+2p_1 \over 2}\right)=
{\beta^{\mu} \alpha^{\mu+1/2} \over 2^{\mu+1} \sqrt{\pi}}\;
\int_0^{\infty}dt\;e^{-(p_1+\alpha/2)t}G(t)\;\;\times
&\nonumber\\&
(t^2+\beta^2)^{-\mu/2-1/4}\; K_{\mu+1/2}\left({\alpha \sqrt{t^2+\beta^2}
\over 2} \right) \quad.
\label{sevennine}&\end{eqnarray}
In Eq.\ (\ref{sevennine}) one can replace $\alpha/2$ by $\alpha$ and
$\alpha/2+p_1$ by $p$, which we do henceforth.
Multiplying both sides of Eq.\ (\ref{sevennine}) by $\beta^{\nu}$ and
integrating over $\beta$ from zero to $\infty$, we get
\begin{eqnarray}&
\int_0^{\infty}ds\;{s^{(\mu-\nu)} \over \sqrt{\alpha^2+s^2}}\;
F\left(\sqrt{\alpha^2+s^2}+p \right)=
{\Gamma((\mu-\nu+1)/2)\; \alpha^{(\mu-\nu)/2} \over 2\sqrt{\pi}}\;\;\times
&\nonumber\\&
2^{(\mu-\nu)/2}\int_0^{\infty}dt\;t^{(\nu-\mu)/2}\;e^{-pt}
G(t) K_{(\mu-\nu)/2}\left(\alpha t \right) \quad,
\label{sevenninea}&\end{eqnarray}
where ${\rm Re}\;(\mu+\nu)>-1$ and ${\rm Re}\;\nu<1/2$.

We now apply the above procedure 
to obtain an apparently new result. Consider the
following fractional integral
\begin{eqnarray}&
I=\int_a^{\infty}dt\;e^{-pt}(t-a)^{\mu-1/2}\;
(t^2+\beta^2)^{-\mu/2}K_{\mu}(p \sqrt{t^2+\beta^2})\quad.
\label{eighttwo}&\end{eqnarray}
Comparing the above with Eq.\ (\ref{sevennine}) we see that
$G(t)=(t-a)^{\mu-1/2}$ and hence, from No. 2.1.2.10 of
Ref.\ \cite{pru92}, $F(x)=\Gamma(\mu+1/2)x^{-(\mu+1/2)}\exp(-ax)$.
Thus, Eq.\ (\ref{sevennine}) yields
\begin{eqnarray}&
I=\sqrt{{\pi \over 2}}\;{e^{-ap}\;\Gamma(\mu+1/2) \over \beta^{\mu-1/2} 
p^{\mu}}\;\int_0^{\infty}ds\;s^{\mu+1/2} J_{\mu-1/2}(\beta s)
\;{e^{-a\sqrt{p^2+s^2}} \over \sqrt{p^2+s^2}}\;=
&\nonumber\\&
\sqrt{{\pi \over 2}}\;{\Gamma(\mu+1/2) \over p^{\mu}}\;
{\beta^{1/2}e^{-p(\sqrt{a^2+\beta^2}+a)} \over (\sqrt{a^2+\beta^2}+a)^{\mu-1/2}
\sqrt{a^2+\beta^2}}\quad,
\label{eightthree}&\end{eqnarray}
where ${\rm Re}\;\mu>-1/2$ and we have used No.\ 2.2.10.13 from Ref.\ \cite{pru86}.

Now consider the following integral where $a>0$ and $b > 0$, 
\begin{eqnarray}&
I=\int_{b}^{\infty}dy\; {e^{-a^2/4(y+b)} \over (y+b)^{3/2}}
\;(y^2-b^2)^{(\mu-\nu-1)/2} \quad.
\label{eightfour}&\end{eqnarray}
If we make the change of variable, $x=(y^2-b^2)^{1/2}$, then
the integral is in a form given by the lhs of Eq.\ (\ref{sevenninea})
where $F(x)=x^{-3/2}\exp(-a^2/4x)$. By utilising the fact that the
inverse LT of $F(x)$ is given as No. 2.2.2.3 in Prudnikov et al
\cite{pru92a}, we find that Eq.\ (\ref{eightfour}) becomes
\begin{eqnarray}&
I={(2b)^{(\mu-\nu)/2} \Gamma((\mu-\nu+1)/2) \over \pi a}\;
\int_0^{\infty}dt\;t^{(\nu-\mu)/2}e^{-b t} \sin(a\sqrt{t})\;
K_{(\mu-\nu)/2}(b t)\;=
&\nonumber\\&
{(2b)^{\mu-\nu-3/2}\Gamma((\mu-\nu+1)/2) \Gamma(\nu-\mu+3/2) \over
2\; \Gamma((\mu-\nu)/2+2)}\;_1F_1\left(\nu-\mu+{3 \over 2};
{\nu-\mu \over 2}+2;-{a^2 \over 8b}\right)\;\;,
\label{eightfive}&\end{eqnarray}
where $| {\rm Re}\;(\mu-\nu)/2|<{\rm Re}\;(\nu-\mu+3)/2$
and we have used No. 2.16.18.2 from Ref.\ \cite{pru86}.

Let us return to the theorem given in  Sec.\ 2 and examine
the denominator of the fractional integral as given by Eq.\ (\ref{two}).
We see that an inversion in the power $\mu$, which we aim to exploit
in the remainder of this paper, occurs when we evaluate
the integral by using Eq.\ (\ref{three}). So far, we have been concerned
with situations where the integrals have been defined but what happens
if we replace $\mu$ by $\mu-n$ in our theorem and let $n$ range from
zero to $\infty$. Clearly, it can be seen that the rhs of Eq.
(\ref{three}) can be become divergent depending, of course, on the
behaviour of $G(t)$. To overcome this potential problem, we can
use the work of Lighthill \cite{lig59} and Ninham \cite{nin66}, who give
interpretations to the class of divergent integral encountered when we
consider a range of values for $n$ after substituting $\mu$ by $\mu-n$ in
our theorem.

If we put $\mu=\mu-n$ in our theorem, then after multiplying both
sides by $(-1)^n \beta^n/n!$ where ${\rm Re}\;\beta>0$, we obtain
\begin{eqnarray}&
{(-1)^n \beta^n \over n! \Gamma(n-\mu+1)}\int_p^{\infty}dx\;
{F(x) \over (x-p)^{\mu-n}}=
{(-1)^n \beta^n \over  n!} \int_0^{\infty}dt\;
t^{\mu-n-1} e^{-pt} G(t) \quad.
\label{eightsix}&\end{eqnarray}
In terms of generalised functions the integral on the rhs of
Eq.\ (\ref{eightsix}) can be written as
\begin{eqnarray}&
I=\int_{-\infty}^{\infty}dt\;t^{\mu-n-1}e^{-pt} H(t)G(t)  \quad,
\label{eightseven}&\end{eqnarray}
where $H(t)$ is the Heaviside step-function. Lighthill \cite{lig59} shows
that if $t^{\mu-1}G(t)$ is a good function, then the above integral yields
the ordinary formula for repeated integration by parts if all the infinite
contributions from the lower limit of the rhs of Eq.\ (\ref{eightsix})
are omitted. As a consequence within the context of generalised functions,
only the finite part to a divergent integral need be kept. This finite part
is referred to as the Hadamard part by Ninham \cite{nin66}, since Hadamard
showed that the finite part obeys many of the ordinary rules of integration.
Ninham's approach to divergent integrals of the form given by Eq.
(\ref{eightseven}) is different from Lighthill's in that he is able to
express such integrals in terms of a sum of function evaluations over an
arbitrary partition of the integration interval, thereby making a connection
with the concept of a Riemann integral as a limit sum. Ultimately, his
approach coincides with Lighthill's but for it to be applicable to integrals
such as Eq.\ (\ref{eightseven}), one needs to make a change of variable such
as $y=p/(t+a)$ in this integral. Then the singularity is shifted to $y=1/a$.
As a numerical example, Ninham chooses the beta function integral 
\begin{eqnarray}&
I(\alpha,\beta)=\int_0^1dt\;t^{\alpha}(1-t)^{\beta}={\Gamma(\alpha+1)
\Gamma(\beta+1) \over \Gamma(\alpha+\beta+2)}\quad.
\label{eighteight}&\end{eqnarray}
According to his prescription, this integral should yield a result
for all $\alpha$ and $\beta$ and for convenience he puts $\alpha=-3/2$
and $\beta=-3/2$. The Hadamard contribution is the final expression
on the rhs of Eq.\ (\ref{eighteight}) and yields zero for these values of
$\alpha$ and $\beta$. Ninham shows that by treating the integral as a sum
over each and every partition of the interval [0,1], denoted by $Rf$, and
an error term associated with the partition, denoted by $Ef$, the integral
virtually yields zero. A slight error arises from the rounding-off of $Rf$
and from the truncation of the asymptotic series for $Ef$.

From the above, we can see that provided it is defined, the lhs of Eq.\
(\ref{eightsix}) represents the Hadamard contribution for the integral
on the rhs. If we sum both sides of Eq.\ (\ref{eightsix}) over
$n$ from zero to $\infty$, then after interchanging the order of the
summations and integrations we get
\begin{eqnarray}&
\int_0^{\infty}dt\;t^{-\mu-1}e^{-pt-\beta/t}\; G(t)=2\beta^{-\mu/2}
\int_0^{\infty}dy\;y^{\mu+1}J_{\mu}(2\sqrt{\beta}\;y)F(y^2+p) \quad,
\label{eightnine}&\end{eqnarray}
where we have put $x=y^2+p$ and ${\rm Re}\;\mu>-1$. For the special case
of $\mu=-1/2$, Eq.\ (\ref{eightnine}) reduces to
\begin{eqnarray}&
\int_0^{\infty}dt\;t^{-1/2}e^{-pt-\beta^2/4t}\;G(t)={2 \over \sqrt{\pi}}
\int_0^{\infty}dy\;F(y^2+p) \cos(\beta y)\quad,
\label{ninezero}&\end{eqnarray}
where we have put $\beta=\beta^2/4$. Note that in the above results our
summation procedure has transformed the singularity at $t=0$ from
either a pole or a branch point, depending upon the value of $\mu$, 
into an essential singularity, which has a less harmful effect on
the evaluation of integrals. We shall see shortly that the summation
of divergent integrals can also lead to a shifting of the singular
point away from the integration interval. This technique of summing 
divergent integrals to obtain a convergent result is the essence of
the renormalisation or regularisation technique used so often in
theoretical physics. In addition, if we had multiplied by $\beta$
instead of $-\beta$, then $\beta$ would have to be replaced by $-\beta$ on
the lhs of Eq.\ (\ref{eightnine}) whilst the Bessel function on the rhs
would have to be replaced by a modified Bessel function. 

To convince the reader of the validity of the above results, we now
consider some examples. Let us put $G(t)$ equal to unity, so that 
the lhs of Eq.\ (\ref{eightnine}) yields $2(\beta/p)^{-\mu/2}K_{\mu}
(2\sqrt{\beta p})$ by No. 3.471.9 of Ref.\ \cite{gra80}. Since $G(t)=1$,
$F(y^2+p)=(y^2+p)^{-1}$ and hence by using No. 6.565.4 of the same
reference, the rhs of Eq.\ (\ref{eightnine}) can be shown to give the
same result as the lhs. On the other hand, if we let $F(x)=x^{-\nu}$
where ${\rm Re}\;\nu>0$, then the lhs of Eq.\ (\ref{ninezero}) yields
\begin{eqnarray}&
\int_0^{\infty}dt\;t^{\nu-3/2}e^{-pt-\beta^2/4t}/\Gamma(\nu)=
{\beta^{\nu-1/2} p^{1/4-\nu/2} \over 2^{\nu-3/2}\Gamma(\nu)}\;
K_{\nu-1/2}(\beta\sqrt{p})\quad.
\label{nineone}&\end{eqnarray}
To obtain this result we have used No. 2.1.1.1 from Ref.\ \cite{pru92a}.
Introducing $(y^2+p)^{-\nu}$ into the rhs of Eq. \ref{ninezero}
gives the answer above after utilising No. 8.432.5 from Gradshteyn
and Ryzhik \cite{gra80}.

A results obtained from Eqs.\ (\ref{eightnine}) and (\ref{ninezero}),
which appears  to be new, is
\begin{eqnarray}&
\int_0^{\infty}dy\;{\cos(\beta y) \over (y^2+p)^{\nu+1/2}}\;
_2F_1\left({\nu \over 2}+{1 \over 4},{\nu \over 2}+{3 \over 4};
\nu+1;-{a^2 \over (y^2+p)^2}\right)= \sqrt{\pi}\;\left({2 \over a}
\right)^{\nu}\;\;\times
&\nonumber\\&
{\Gamma(\nu+1) \over \Gamma(\nu+1/2)}\;J_{\nu}\left(\beta 
\left(\sqrt{p^2+a^2}-p\right)^{1/2}\right) K_{\nu}\left(\beta
\left(\sqrt{p^2+a^2}+p\right)^{1/2}\right)
&\nonumber\\&
[{\rm Re}\;p,\;{\rm Re}\;\beta^2 >0,\;{\rm Re}\;\nu>-1/2]\quad,
\label{ninethree}&\end{eqnarray}

To study a situation where the summation of divergent integrals
produces a singularity shifted away from the integration interval
we multiply Eq.\ (\ref{eightsix}) by $\Gamma(n-\mu+1)$, replace
$\mu$ by $-\mu$ and then sum over $n$ from zero to $\infty$.
Then after interchanging the summations and integrations we get
\begin{eqnarray}&
\int_p^{\infty}dx\;(x-p)^{\mu}F(x)e^{-\beta(x-p)}= \Gamma(1+\mu)
\int_0^{\infty}dt\;e^{-pt}G(t) t^{-\mu-1}\;_1F_0(1+\mu;-\beta/t)\quad,
\label{ninefive}&\end{eqnarray}
where ${\rm Re}\;\mu>-1$. Introducing No. 7.3.1.1 from Prudnikov et al
\cite{pru90} into Eq.\ (\ref{ninefive}) yields
\begin{eqnarray}&
\int_p^{\infty}dx\;(x-p)^{\mu}F(x)e^{-\beta(x-p)}= \Gamma(1+\mu)
\int_0^{\infty}dt\;{e^{-pt}\; G(t) \over (t+\beta)^{1+\mu}}\quad.
\label{ninesix}&\end{eqnarray} 
Hence, we can see that the singularity has been shifted to $-\beta$
and is now outside the region of integration. Eq.\ (\ref{ninesix})
can be verified by letting $F(x)=x^{-1}$ or $G(t)=1$. To show that the
lhs is indeed equal to the rhs, put $y=x-p$ to get
\begin{eqnarray}&
\int_0^{\infty}dy\;{y^{\mu}e^{-\beta y} \over y+p}=\int_0^{\infty}dt\;
e^{-pt}\int_0^{\infty}dy\;y^{\mu}e^{-(\beta+t)y}\quad,
\label{nineseven}&\end{eqnarray}
which is just the rhs of Eq.\ (\ref{ninesix}). By multiplying Eq.\
(\ref{eightsix}) by $\Gamma(n+\nu+1)$, we can generalise Eq.\ (\ref{ninefive})
to
\begin{eqnarray}&
\int_0^{\infty}dy\;y^{\mu}F(y+p)\;_1F_1\left(\nu+1;\mu+1;-\beta y\right)=
\Gamma(\mu+1) \int_0^{\infty}dt\; G(t)\;{t^{\nu-\mu}e^{-pt} \over
(t+\beta)^{\nu+1}}\;\;,
\label{ninesevena}&\end{eqnarray}
where ${\rm Re}\;\mu>-1$. In addition, if we put $\mu=n+\mu$ and multiply
both sides of Eq.\ (\ref{ninesix}) by $(-1)^n \gamma^{2n+\mu}/2^{2n+\mu}n!
\Gamma(n+\mu+1)$, then after interchanging the summations and integrations
we obtain
\begin{eqnarray}&
\int_0^{\infty}ds\;s^{\mu+1}J_{\mu}(\gamma s)e^{-\beta s^2} F(s^2+p)=
{\gamma^{\mu} \over 2^{\mu+1}}\int_0^{\infty}dt\;G(t)\;{e^{-pt-
\gamma^2/4(t+\beta)} \over (t+\beta)^{\mu+1}} \quad.
\label{nineeight}&\end{eqnarray}

We can still go further with Eq.\ (\ref{eightsix}). By dividing both
sides by $\Gamma(n+\nu+1)$ and then carrying out our summation procedure
we find
\begin{eqnarray}&
\int_p^{\infty}dx\;(x-p)^{\mu} F(x)\;_0F_2\left(\mu+1,\nu+1;
-\beta(x-p)\right)= \beta^{-\nu/2}\;\Gamma(\mu+1)\; \Gamma(\nu+1)\;\times
 &\nonumber\\&
\int_0^{\infty}dt\;t^{\mu-\nu/2-1}\; e^{-p/t}\;G(t^{-1})\;
J_\nu\left(2\sqrt{\beta t}\right)\quad,
\label{onezerozero}&\end{eqnarray}
where ${\rm Re}\;\mu>-1$. By putting $F(x)=x^{-1}$ and using No.\ 2.12.9.14
from Ref.\ \cite{pru86} with the above result we get
\begin{eqnarray}&
\int_0^{\infty}dy\; {y^{\mu} \over (y+p)}\;_0F_2\left(1+\mu,1+\nu;
-\beta y\right)= \Gamma(\nu+1) \Gamma(\mu+1)\;
\left[ \beta^{-\mu}\;{\Gamma(\mu) \over \Gamma(1+\nu-\mu)}\;\;\times\right.
&\nonumber\\&
\left. _0F_2\left(1-\mu,1+\nu-\mu;\beta p\right)+ p^{\mu}\; {\Gamma(-\mu) \over
\Gamma(\nu+1)}\;_0F_2\left(\mu+1,\nu+1;\beta p\right)\right] \quad,
\label{onezeroone}&\end{eqnarray}
where ${\rm Re}\;(2\mu-\nu)<3/2$, ${\rm Re}\;p>0$ and
$\sqrt{\beta}>0$. 

To complete this section we utilise the preceding material to obtain 
the following general result for ${\rm Re}\;\mu>-1$
\begin{eqnarray}&
\int_p^{\infty}dx\;(x-p)^{\mu}F(x)\;_pF_{q+1}\left(\alpha_1,...,\alpha_p;
\mu+1,\gamma_1,...,\gamma_{q};-\beta(x-p)\right)=\Gamma(\mu+1)\;\;\times
&\nonumber\\&
\int_0^{\infty}dt\;t^{-\mu-1}e^{-pt}G(t)\;_pF_q\left(\alpha_1,...,\alpha_p;
\gamma_1,...,\gamma_q;-{\beta \over t}\right)\quad.
\label{onezerotwo}&\end{eqnarray}
This result represents the generalisation of the theorem
presented in  Sec.\ 2 and is not always restricted to ${\rm Re}\;\beta>0$. 
Using the tabulations of hypergeometric functions given in Ref.
\cite{pru90}, we obtain the following results from the above equation
\begin{eqnarray}&
\int_0^{\infty}dy\;y^{3/2-b}\; F(y+p) J_{b-1}\left(\beta \sqrt{y}\right)
J_{2-b}\left(\beta \sqrt{y}\right)=\beta^{-1}\;{b \over 3 \Gamma(b)}\;
\int_0^{\infty}dt\;t^{b-2}\;e^{-pt}\;\times
&\nonumber\\&
G(t)\left(1-\; _1F_1\left({3 \over 2};b;-{\beta^2 \over t}\right)\right)
&\nonumber\\&
[{\rm Re}\;b >0]\quad,
\label{onezerothree}&\end{eqnarray}
and
\begin{eqnarray}&
\int_p^{\infty}dx\;F(x) J^2_{a-1/2}\left(\sqrt{\beta (x-p)}\right)=
\int_0^{\infty}dt\;t^{-2}\;e^{-pt-\beta/2t}\;
G(t)\;I_{a-1/2}\left(-{\beta \over 2t}\right)
&\nonumber\\&
[{\rm Re}\;a>-1/2]\quad.
\label{onezerofive}&\end{eqnarray}
Putting $b=3/2$ in Eq.\ (\ref{onezerothree}) eventually yields the
result given by Eq.\ (\ref{ninezero}).

\section{Conclusion}
In this paper we have presented a theorem/technique derived from Laplace
transform theory that facilitates the evaluation of Weyl fractional integrals
by transforming them into known integrals. As a consequence, we 
were able to present many new fractional integrals not previously
evaluated in the standard tables of integrals, Refs. \cite{gra80}-
\cite{ape83}, \cite{pru86a} and \cite{pru90}. Some of these results
appear in Sec.\ 2 where we introduced and proved the theorem.
Further examples are presented in Appendix A.
In addition, the technique is particularly useful for the
numerical evaluation of fractional integrals since it is able to transform
slowly converging integrals into rapidly converging ones.

We showed in Sec.\ 3 how the theorem could be extended to integrals
not of a fractional form but concentrated upon integrals which yielded
decaying exponential-like behaviour after transformation. As a result,
we were able to evaluate many new integrals, which are listed in Appendix B.
In  Sec.\ 4 we considered several extensions of our theorem.
During the course of our study we 
found that by utilising Feynman's integral from quantum electrodynamics 
together with the techniques in this paper we could derive an
interesting
class of results, anumber of which are given 
in Appendix C. Finally, by using the divergent integral theory of Lighthill
\cite{lig59} and Ninham \cite{nin66}, we were able to show that our
theorem was only a special case of the result given by Eq.\ (\ref{onezerotwo}).

\section*{aknowledgements}
The authors are grateful to the Australian Research Council for
supporting this work. One of us (V.K.) thanks Professor B.W. Ninham
for bringing Ref.\ \cite{nin66} to his attention during a visit
to the Department of Applied Mathematics at the Institute of Advanced
Studies, Australian National University, Canberra.

\section{Appendix A}

We begin this appendix by utilising the theorem in  Sec.\ 2 to
evaluate a few more fractional integrals. In listing these results
we also indicate which references were used, so that should
any errors appear, the reader may trace their origin. 

\vspace{0.5cm}
\begin{eqnarray}&
\int_p^{\infty}dx\;(\sqrt{x}+a)^{\nu}(x-p)^{\mu}={-\nu (2a)^{2\mu+\nu+2}
\;\Gamma(\mu+1) \Gamma(-2\mu-\nu-2) \over \Gamma(-\mu-\nu)}\;\;\times
&\nonumber\\&
_2F_1\left(-\mu-{\nu \over 2}-1,-\mu-{\nu+1 \over 2};
-\mu-\nu;1-{p \over a^{2}}\right)
&\nonumber\\&
[{\rm Re}\;\mu >-1,\;{\rm Re}\;\nu <0,\;{\rm Re}\;\mu<- {\rm Re}\;
\nu/2-1,\; |{\rm arg}\;a|<3 \pi/4] 
\label{athree}&\end{eqnarray}
[No. 2.1.7.23 Ref.\ \cite{pru92a}, No. 2.11.3.2 Ref.\ \cite{pru86}],

\vspace{0.5cm}
\begin{eqnarray}&
\int_p^{\infty}dx\;\left({(\sqrt{x^2+a^2} \pm x)^{\nu} \over
\sqrt{x^2+a^2}}\right)(x-p)^{\mu}={a^{\nu}\;\Gamma(\mu+1)\;
\Gamma(\mp \nu-\mu) \over p^{-\mu} \Gamma(\mp \nu+1)}\;\;\times
&\nonumber\\&
\left({a \over 2p}\right)^{\mp \nu}\;_2F_1\left({-\mu \mp \nu \over 2},
{-\mu \mp \nu+1 \over 2};\mp \nu+1;-{a^2 \over p^2}\right)
&\nonumber\\&
[{\rm Re}\;\mu>-1,\;\mp {\rm Re}\;\nu >{\rm Re}\;\mu,\; p>|{\rm Im}\;a|] 
\label{afive}&\end{eqnarray}
[No. 2.1.9.15 Ref.\ \cite{pru92a}, No. 2.12.8.4 Ref.\ \cite{pru86}],

\vspace{0.5cm}
\begin{eqnarray}&
\int_p^{\infty}dx\;{(x-p)^{\mu} \over \sqrt{x^2-a^2}}\;\left(
(x+\sqrt{x^2-a^2})^{\nu}-(x-\sqrt{x^2-a^2})^{\nu}\right)=
{2^{\nu+1}\;a^{2\nu} \sin(\nu \pi) \over \sqrt{\pi}\;(p+a)^{\nu-\mu}}
\;\;\times
&\nonumber\\&
\Gamma(\mu+1)\Gamma(\nu-\mu)\;{\Gamma(-\nu-\mu) \over \Gamma(1/2-\mu)}\;
_2F_1\left(\nu-\mu,\nu+{1 \over 2};{1 \over 2}-\mu;
{p-a \over p+a}\right)
&\nonumber\\&
[| {\rm Re}\;\nu|<1,\;-{\rm Re}\;\mu> |{\rm Re}\;\nu|,\;
{\rm Re}\;(a+p)>0] 
\label{asix}&\end{eqnarray}
[No. 2.1.9.24 Ref.\ \cite{pru92a}, No. 2.16.6.3 Ref.\ \cite{pru86}],

\vspace{0.5cm}
\begin{eqnarray}&
\int_p^{\infty}dx\;(x-p)^{\mu}\;e^{-ax^{1/3}}={\Gamma(\mu+1) \over \pi}
\left[{3^{3\mu+8} \over 2\; a^{3\mu+15/2}}\;\Gamma\left(\mu+{5 \over 3}
\right)\Gamma\left(\mu+{4 \over 3}\right)\right.\;\;\times
&\nonumber\\&
_0F_2\left(-\mu-{2 \over 3},-\mu-{1 \over 3};
-{pa^3 \over 27}\right)+ {3^3p^{\mu+5/3} \over 2 a^{5/2}}\;
\Gamma\left(-{1 \over 3}\right)\Gamma\left(-\mu-{5 \over 3}\right)\;\times
&\nonumber\\&
\left. _0F_2\left(\mu+{8 \over 3}, {4 \over 3};-{pa^3 \over 27}\right)+
{3^4 p^{\mu+4/3} \over 2 a^{7/2}}\;\Gamma\left({1 \over 3}\right)
\Gamma\left(-\mu-{4 \over 3}\right)\;_0F_2\left(\mu+{7 \over 3},
{2 \over 3};-{pa^3 \over 27}\right)\right] 
&\nonumber\\&
[|{\rm arg}\;a|< \pi/3,\;{\rm Re}\;\mu>-1] 
\label{anine}&\end{eqnarray}
[No. 2.2.1.4 Ref.\ \cite{pru92a}, No. 2.16.8.13 Ref.\ \cite{pru86}],

\vspace{0.5cm}
\begin{eqnarray}&
\int_p^{\infty}dx\;{(x-p)^{\mu} \over x}\;e^{-ax^{1/2}}=
{\Gamma(\mu+1) \over 2\sqrt{\pi}}\left[ap^{\mu+1/2}\Gamma(-(\mu+1/2))
\;\;\times\right.
&\nonumber\\&
\left. _1F_2\left({1 \over 2};{3 \over 2},\mu+{3 \over 2};
{a^2 p \over 4}\right)-{2^{2\mu}\Gamma(\mu+1/2) \over  \mu a^{2\mu}}
\;_1F_2\left(-\mu;1-\mu,{1 \over 2}-\mu;{a^2 p \over 4}\right)\right]
&\nonumber\\&
[{\rm Re}\;a^2>0,\; -1<{\rm Re}\;\mu<0]
\label{aonezero}&\end{eqnarray}
[No. 2.2.1.16 Ref.\ \cite{pru92a}, No. 2.8.5.15 Ref.\ \cite{pru86}],

\vspace{0.5cm}
\begin{eqnarray}&
\int_p^{\infty}dx\;(x-p)^{\mu}\;e^{a^2x^2}{\rm erfc}(ax)=
{\Gamma(\mu+1)\Gamma(-\mu) \over 2^{\mu/2}\;a^{\mu+1} \;\sqrt{\pi}}\;
e^{a^2p^2/2}\;D_{\mu}\left(\sqrt{2}ap\right)
&\nonumber\\&
[-1< {\rm Re}\;\mu<0,\;{\rm Re}\;a^2>0]
\label{aonesix}&\end{eqnarray}
[No. 3.7.3.2 Ref.\ \cite{pru92a}, No. 2.3.15.3 Ref.\ \cite{pru86a}],  

\vspace{0.5cm}
\begin{eqnarray}&
\int_p^{\infty}dx\;{(x-p)^{\mu} \over x^{2a}}(\omega^2+x^2)^{2a-c+1/2}
\;_2F_1(a,a+1/2;c;-\omega^2/x^2)={\Gamma(\mu+1)\over \Gamma(2c-2a-1)}
\;\;\times
&\nonumber\\&
{\Gamma(2c-2a-\mu-2) \over  p^{2c-2a-\mu-2}}\;
_2F_1(c-a-\mu /2-1,c-a-\mu /2-1/2;c;-\omega^2/p^2)
&\nonumber\\&
[{\rm Re}\;c>{\rm Max}[{\rm Re}(\mu/2+a)+1,\;{\rm Re}\;a+1/2],\;
p>|{\rm Im}\;\omega|]
\label{aoneeight}&\end{eqnarray}
[No. 3.35.1.22 Ref.\ \cite{pru92a}, No. 2.12.8.4 Ref.\ \cite{pru86}],

\vspace{0.5cm}
\begin{eqnarray}&
\int_p^{\infty}dx\;{(x-p)^{\mu} \over x^{2a}}\;_2F_1(a,a+1/2;c;
-\omega/x^2)={\Gamma(\mu+1) \over \Gamma(2a)}\;\;\times
&\nonumber\\&
{\Gamma(2a-\mu-1) \over p^{2a-\mu-1}}\;_2F_1\left({2a-\mu-1 \over 2},
{2a+\mu-c \over 2};c;-{\omega \over p^2}\right)
&\nonumber\\&
[{\rm Re}\;a>0,\;p>|{\rm Im}\;\sqrt{\omega}|,\;{\rm Re}\;
(2a-\mu-1)>0,\;{\rm Re}\;\mu>-1]
\label{aonenine}&\end{eqnarray}
[No. 3.35.1.17 Ref.\ \cite{pru92a}, No. 2.12.8.4 Ref.\ \cite{pru86}],

\vspace{0.5cm}
\begin{eqnarray}&
\int_0^1dx\;x^{\gamma-1/2}(1-x)^{-1/2}\;_{r+n}F_s\left(\alpha_1,
...,\alpha_r,{\gamma+1 \over n},...,{\gamma+n \over n};\beta_1,
...,\beta_s;-\omega (nx)^n\right)=
&\nonumber\\&
\sqrt{\pi}\;{\Gamma(\gamma+1/2) \over \Gamma(\gamma+1)}\;
_{r+n}F_s\left(\alpha_1,...,\alpha_r,{\gamma+1/2 \over n},...,
{\gamma+n-1/2 \over n};\beta_1,...,\beta_s;-\omega n^n\right)
&\nonumber\\&
[r+n \leq s+1,\;{\rm Re}\;\gamma >1/2,\;\alpha_i \neq 0,-1,-2,...]
\label{aoneninea}&\end{eqnarray}
[No. 3.36.1.11 Ref.\ \cite{pru92a}, No. 21.26 Ref.\ \cite{obe73}],

We also mentioned in  Sec.\ 2 that our theorem can be utilised in the
construction of fractional integrals when it cannot be used to solve them.
As an example we cited Eq.\ (\ref{oneone}), whose more general form, given by 
Eq.(9), was evaluated by using our theorem. We were able to use 
this result to develop some fractional integrals. We should add, however, 
that more fractional integrals can be developed from this result, which we
present here. Thus, by putting $\nu$ equal to $n+\nu$ in Eq.\ (\ref{oneone}),
multiplying both sides by $\beta^n/n!$ and summing from $n=0$ to $\infty$,
we get after interchanging the order of the summation and integration
\begin{eqnarray}&
\int_p^{\infty}{(x+a)^{-\nu} \over (x-p)^{\mu}}\;e^{\beta/(x+a)}=
{\Gamma(1-\mu) \Gamma(\mu+\nu-1) \over (p+a)^{\mu+\nu-1} \Gamma(\nu)}\;
_1F_1\left(\mu+\nu-1;\nu;{\beta \over p+a}\right)\quad,
\label{atwozero}&\end{eqnarray}
where the same conditions apply to $\mu$ and $\nu$ as in Eq.(9).
Similarly, one finds 
\begin{eqnarray}&
\int_p^{\infty}dx\;{(x+a)^{-\nu} \over (x-p)^{\mu}}\;
\exp\left({\beta (x-p) \over x+a}\right)= {\Gamma(1-\mu)
\Gamma(\mu+\nu-1) \over \Gamma(\nu)\;(p+a)^{\mu+\nu-1}}\;
_1F_1(1-\mu;\nu;\beta) \quad,
\label{atwoone}&\end{eqnarray}
\begin{eqnarray}&
\int_p^{\infty}dx\;{(x+a)^{(\mu-\nu)/2-1} \over (x-p)^{\mu}}\;
J_{(\mu+\nu)/2-1}\left({2\sqrt{\beta} \over x+a}\right)=
{2^{\mu-1} \beta^{(\mu+\nu)/4-1/2} \over (p+a)^{\mu+\nu-1}}\;
{\Gamma(1-\mu) \over \Gamma(\nu)}\;\;\times
&\nonumber\\&
{\Gamma((\mu+\nu-1)/2) \over \Gamma((\nu+1)/2)}\;
_1F_2\left({\mu+\nu-1 \over 2};{\nu \over 2},{\nu+1 \over 2};
{-\beta \over (p+a)^2}\right)\quad,
\label{atwotwo}&\end{eqnarray}
\begin{eqnarray}&
\int_p^{\infty}dx\;{(x+a)^{-\nu} \over (x-p)^{\mu}}\;
\exp\left(-\beta (x-p)^2/(x+a)^2 \right)={2^{1-\mu-\nu} \Gamma((1-\mu)/2)
\over (p+a)^{\mu+\nu+1} \Gamma((\nu+1)/2)}\;\;\times
&\nonumber\\&
{\Gamma((1-\mu)/2) \over \Gamma(\nu/2)}\;_2F_2\left({1 \over 2}-{\mu
\over 2},1-{\mu \over 2};{\nu \over 2},{\nu+1 \over 2};-\beta\right)\quad,
\label{atwothree}&\end{eqnarray}
\begin{eqnarray}&
\int_p^{\infty}dx\;{J_{\nu-1/2}(\beta/(x+a)) \over (x-p)^{1/2}
(x+a)^{3/2}}={\pi \over \sqrt{2}\;(p+a)}\;J_{\nu/2}\left(
{\beta \over 2(p+a)}\right) J_{(\nu-1)/2}\left({\beta \over
2(p+a)}\right) \;,
\label{atwofour}&\end{eqnarray}
\begin{eqnarray}&
\int_p^{\infty}dx\;{(x+a)^{-\nu} \over (x-p)^{\mu}}\;
_0F_3\left({\nu+\mu-1 \over 2},{\nu+\mu \over 2},\nu;{\beta^4 \over
4(x+a)^2}\right)={\Gamma(1-\mu) \over \beta^{2\nu-2}}\;\times
&\nonumber\\&
{\Gamma(\nu+\mu-1) \Gamma(\nu) \over \sqrt{\pi}\;(p+a)^{\mu}}
\left(ber_{\nu-1}^2\left({\beta \over 2 \sqrt{p+a}}\right) +
bei_{\nu-1}^2\left({\beta \over 2 \sqrt{p+a}}\right)\right) \quad,
\label{atwofive}&\end{eqnarray}
and
\begin{eqnarray}&
\int_p^{\infty}dx\;{(x+a)^{\mu-1} \over (x-p)^{\mu}}\left[
ber^2_{\nu+\mu-1}\left({\beta \over \sqrt{x+a}}\right)+
bei^2_{\nu+\mu-1}\left({\beta \over \sqrt{x+a}}\right)\right]=
&\nonumber\\&
\left({\beta^2 \over 4(p+a)}\right)^{\nu+\mu-1}\;
{\Gamma(1-\mu) \over \Gamma(\nu) \Gamma(\nu+\mu-1)}
\;_0F_3\left({\nu \over 2},{\nu+1 \over 2},
\nu+\mu-1;{\beta^4 \over 64(p+a)^2}\right)\;.
\label{atwosix}&\end{eqnarray}

We also stated in  Sec.\ 2 that more results could be determined from
the integral given in Eq.\ (\ref{twotwo}) by using the tabulated 
results for $_3F_2$ hypergeometric functions in Ref.\ \cite{pru90}.
For example, if we put $\alpha=-2$ in Eq.\ (\ref{twotwo}), then we get
\begin{eqnarray}&
\int_p^{\infty}dx\;{x^2(x^2+a^2)^{-\nu} \over (x-p)^{\mu}}=
p^{3-\mu-2\nu}\;B(1-\mu,\mu+2\nu-3)\;\;\;\times
&\nonumber\\&
\left[_2F_1\left(\nu+{\mu-3 \over 2},\nu+{\mu \over 2}-1;
\nu-{1 \over 2};-{a^2 \over p^2}\right)-
{(\nu+(\mu-3)/2) (\nu+\mu/2-1) a^2 \over (\nu-1) (\nu-1/2) p^2}\;\;\times\right.
&\nonumber\\&
\left. _2F_1\left(\nu+{\mu-1 \over 2},\nu+{\mu \over 2};\nu+{1 \over 2};-
{a^2 \over p^2}\right)\right] \quad.
\label{atwoseven}&\end{eqnarray} 
If we put $\alpha=1$ and $\mu=1/2-\nu$, then Eq.\ (\ref{twotwo}) yields
\begin{eqnarray}&
\int_p^{\infty}dx\;{x^{-1}(x^2+a^2)^{-\nu} \over (x-p)^{1/2-\nu}}=
p^{-\nu-1/2}\;B(\nu+1/2,\nu+1/2)\;
\left({1+\sqrt{1+a^2/p^2} \over 2}\right)^{-\nu}\;\times
&\nonumber\\&
_2F_1\left(\nu,{1 \over 2};\nu+1;{1+\sqrt{1+a^2/p^2} \over 2}\right)
\quad.
\label{atwoeight}&\end{eqnarray}
Eq.\ ({\ref{atwoeight}) is based on using No.\ 7.4.1.13 in Prudnikov et
al \cite{pru90}, which has another representation given immediately above
it in the same table. Hence, the above result can be expressed in a different
form. On the other hand, if we put $\nu=1$ and $\alpha=2$, then we find
\begin{eqnarray}&
\int_p^{\infty}dx\;{x^{-2} (x^2+a^2)^{-1} \over (x-p)^{\mu}}=
-{6p^2 \over (\mu^2+3\mu+2)a^2}\left[_2F_1\left({\mu+1 \over 2},
{\mu \over 2}+1;{3 \over 2};-{a^2 \over p^2}\right)-1\right] \quad.
\label{atwonine}&\end{eqnarray}

To complete this appendix we now give further results that can be
derived as a result of the evaluation of Eq.\ (\ref{threefive}),
whose answer is expressed as a $_4F_3$ hypergeometric function in
Eq.\ (\ref{threeseven}). Thus, using the results in  Sec.\ 7.5.1 of
Ref.\ \cite{pru90}, one obtains
\begin{eqnarray}&
\int_0^1dy\;{y^{\alpha-1} \over (1-y)^{\mu}}\;_2F_1\left(a,a+{1 \over 2};
{3 \over 2};-{\omega y^2 \over p^2}\right)={-i\;\omega^{-1/2}p(\alpha-\mu)
\Gamma(\alpha) \over 2(2a-1)(\alpha-1)\Gamma(\alpha+1-\mu)}\;\;\times
&\nonumber\\&
\left[_2F_1\left(2a-1,\alpha-1;\alpha-\mu;{i \omega^{1/2} \over p}\right)
-\;_2F_1\left(2a-1,\alpha-1;\alpha-\mu;-{i \omega^{1/2} \over p}\right)
\right]\quad ,
\label{athreezero}&\end{eqnarray}
\begin{eqnarray}&
\int_0^1dy\;{y^{\alpha-1} \over (1-y)^{\mu}}\;_2F_1\left({\alpha \over 2},
1-\mu+{\alpha \over 2};{\alpha +1 \over 2}-\mu;-{\omega y^2 \over p^2}\right)
={\Gamma(\alpha) (1-z)^{\alpha/2} \over \Gamma(\alpha+1-\mu)}
\;\;\times
&\nonumber\\&
_2F_1\left({\alpha \over 2},{\alpha +1 \over 2}-\mu;
\alpha+1-\mu;z\right)\;_2F_1\left({\alpha \over 2},{\alpha +1 \over 2};
\alpha+1-\mu;z\right) \;\;,
\label{athreeone}&\end{eqnarray}
and
\begin{eqnarray}&
\int_0^1dy\;{1-y \over y^{3/2}}\;\arctan\left({\omega^{1/2}y \over p}
\right)={3 \omega^{1/2} \over 4 p}\;z^{-3/4}\biggl[(1+z^{1/4})^2
\ln(1+z^{1/4})\;\;-
&\nonumber\\&
(1-z^{1/4})^2 \ln(1-z^{1/4})\biggr]-{3 \omega^{1/2} \over 2 p \sqrt{z}}\;
\ln(1+\sqrt{z})+{3 \omega^{1/2} \over 2 p z^{3/4}} \;(\sqrt{z}-1)
\arctan(z^{1/4})\quad.
\label{athreetwo}&\end{eqnarray}
In Eq.\ (\ref{athreeone}) $z$ is a solution of $z^2=-4\omega(z-1)/p^2$
while in Eq.\ (\ref{athreetwo}) $z=\exp(i\pi)\omega/p^2$.

\section{Appendix B}
In this appendix we develop further results based on the material
presented in  Sec.\ 3. We shall concentrate on integrals where part of the
integrand has an inverse LT with some form of exponential behaviour.
As in the previous appendix the list of results presented here is by
no means exhaustive.

We begin with integrals containing the factor of $x^{\nu}\exp(-a/x)$ in
their integrands. Since the LT of this factor is given as No. 2.2.2.1
in Ref.\ \cite{pru92}, we find 
\begin{eqnarray}&
\int_0^{\infty}dx\;x^{\nu}e^{-a/x} F(x)=2a^{(\nu+1)/2}\int_0^{\infty}dt\;
t^{-(\nu+1)/2}\;K_{\nu+1}(2\sqrt{at})\;G(t)\quad,
\label{bone}&\end{eqnarray}
where $G(t)$ is the inverse LT of $F(x)$ and ${\rm Re}\;a >0$. We now
present some new results arising from Eq.\ (\ref{bone}). As before we
list the references used in obtaining these results, so that the reader 
may verify their correctness. Thus we obtain
\begin{eqnarray}&
\int_0^{\infty}dx\;x^{\nu}e^{-a/x-px^{1/2}}= \left[-p\;a^{\nu+3/2}\; 
\Gamma\left(-\nu-{3 \over 2}\right)\;
_0F_2\left({3 \over 2},\nu+{5 \over 2};-ap^2 \right)\;+ \right.
&\nonumber\\&
4 a^{\nu+1} \; \Gamma(-\nu-1)\;
_0F_2\left({1 \over 2},\nu+2;-ap^2\right)+\;
{2^{2\nu+2} \over \sqrt{\pi}\;p^{2\nu+2}}\; \Gamma(\nu+1)\;\;\times
&\nonumber\\&
\Gamma\left(\nu+{3 \over 2}\right)
\left._0F_2\left(-\nu-{1 \over 2}, -\nu;-ap^2\right)\right]
&\nonumber\\&
[{\rm Re}\;p^2 >0,\;{\rm Re}\;\sqrt{a}>0]
\label{btwo}&\end{eqnarray}

[No. 2.2.1.9 Ref.\ \cite{pru92a}, No. 2.16.8.13 Ref.\ \cite{pru86}],

By using No. 2.16.5.4 from Ref.[2], we get the new Laplace transform
\vspace{0.5cm}
\begin{eqnarray}&
\int_0^{\infty}dx\;x^{\lambda-1}e^{-ax}(\sqrt{x^2+1}+1)^{-\mu}=
{\mu \over 4}\left[{\Gamma(\lambda-\nu) \Gamma(\mu/2) a^{\mu-\lambda}
\over  \Gamma((\mu-\nu+1)/2)}\;\;\times\right.
&\nonumber\\&
_2F_3\left({-\lambda \over 2},{\mu \over 2};{1 \over 2},{\mu-\lambda+1 \over 2}+1,
{\mu-\lambda+2 \over 2};-{a^2 \over 4 }\right)\;\;-
{2 \Gamma(\lambda-\mu-1) \Gamma((\mu+1)/2) a^{\nu+2} \over 
\Gamma((\lambda+1)/2)}\;\times
&\nonumber\\&
_2F_3\left(1-{\lambda+1 \over 2}, {\mu+1 \over2};
{3 \over 2},{\mu-\lambda+2 \over 2},{3+\mu-\lambda \over 2};-{a^2 \over 4}\right)
+{\Gamma(\mu-\lambda) \Gamma((\mu-\nu-1)/2) \over 2^{\mu-\lambda}\;
 \Gamma(\mu/2+1)}\;\times
&\nonumber\\&
_2F_3\left(-{\mu \over 2}, {\lambda \over 2};{1 \over 2},{1-\mu+\lambda \over 2},
{2+\lambda-\mu \over 2}; -{a^2 \over 4 }\right)\;-
{a \Gamma(\mu-\lambda+1) \Gamma((\lambda+1)/2) \over 2^{\mu-\lambda-1}
\Gamma((\mu+1)/2)}\;\times
&\nonumber\\&
\left._2F_3\left({1-\mu \over 2}, {\lambda+1 \over 2};{3 \over 2},
{2+\lambda-\mu \over 2},{3-\mu+\lambda \over 2};- {a^2 \over 4 }\right)\right]
&\nonumber\\&
[{\rm Re}\;\mu>0,\;{\rm Im}\;a=0,\;{\rm Re}\;\sqrt{a}>0,\;
\;{\rm Re}\;(\mu+\lambda)>|{\rm Re}\;(\mu+\lambda|)]
\label{bfour}&\end{eqnarray}

[No. 2.1.9.15 Ref.\ \cite{pru92a}, No. 2.16.22.2 Ref.\ \cite{pru86}].

Integrals containing $x^{\nu}(x-a)^{\nu}$ where $a>0$ can be transformed
as follows
\begin{eqnarray}&
\int_a^{\infty}dx\;x^{\nu}(x-a)^{\nu} F(x) ={\Gamma(\nu+1) a^{\nu+1/2}
\over \sqrt{\pi}}\;\int_0^{\infty}dt\;t^{-\nu-1/2}e^{-at/2}
K_{\nu+1/2}(at/2)G(t)\quad,
\label{bfive}&\end{eqnarray}
where ${\rm Re}\;\nu>-1$. Some examples arising from this result are
\begin{eqnarray}&
\int_a^{\infty}dx\;x^{\nu}(x-a)^{\nu}\;{(\sqrt{x^2+b^2}-x)^{1/2} \over
\sqrt{x^2+b^2}}=2^{-1/2} a^{2\nu-1/2}b\;{\Gamma(\nu+1) \Gamma(1/2-2\nu)
\over \Gamma(3/2-\nu)}\;\;\times
&\nonumber\\&
_4F_3\left({3 \over 4},{5 \over 4},{1 \over 4}-\nu,{3 \over 4}-\nu;
{3 \over 2},{3 \over 4}-{\nu \over 2},{5 \over 4}-{\nu \over 2};
-{b^2 \over a^2}\right)
&\nonumber\\&
[b>0,\;-1<{\rm Re}\;\nu<1/4]
\label{bsix}&\end{eqnarray}

[No. 2.1.8.25 Ref.\ \cite{pru92a}, No. 2.16.18.1 Ref.\ \cite{pru86}],

\vspace{0.5cm}
\begin{eqnarray}&
\int_a^{\infty}dx\;x^{\nu}(x-a)^{\nu}\exp(-bx^{1/2})={\Gamma(\nu+1)
\over \pi}\;\left[{\sqrt{\pi}\;a^{2\nu+1} \over 2^{2\nu+2}}\;
\Gamma\left(-\nu-{1 \over 2}\right)\;\times \right.
&\nonumber\\&
_1F_2\left(\nu+1;2\nu+2,{1 \over 2};{ab^2 \over 4}\right) + {2^{6\nu+2}
\over b^{2\nu+2}}\; \Gamma\left(\nu +{1 \over 2}\right) \Gamma\left(2\nu+
{3 \over 2}\right)\;\;\times
&\nonumber\\&
_1F_2\left(-\nu;-2\nu,-2\nu-{1 \over 2};{ab^2 \over
4}\right)+
\pi a^{2\nu+3/2}b\; {\Gamma(-2\nu-3/2)
\over \Gamma(-\nu-1/2)}\;\times
&\nonumber\\&
\left._1F_2\left(\nu+{3 \over 2};2\nu+{5 \over 2},
{3 \over 2};{ab^2 \over 4}\right)\right]
&\nonumber\\&
[{\rm Re}\;b^2 >0,\;{\rm Re}\;\nu>-1]
\label{bseven}&\end{eqnarray}

[No. 2.2.1.9 Ref.\ \cite{pru92a}, No. 2.16.9.3 Ref.\ \cite{pru86}],
\newline

Integrals of the form
\begin{eqnarray}&
I=\int_a^{\infty}dx\;{F(x) \over \sqrt{x^2-a^2}}\;\left[
(x+\sqrt{x^2-a^2})^{\nu}+(x-\sqrt{x^2-a^2})^{\nu}\right]
\label{bnine}&\end{eqnarray}
can be written as
\begin{eqnarray}&
I=2a^{\nu}\int_0^{\infty}dt\;G(t) K_{\nu}(at) \quad.
\label{bten}&\end{eqnarray}
Some interesting examples from the above result are
\begin{eqnarray}&
\int_a^{\infty}dx\;{(\sqrt{x^2+b^2}-x)^{\mu} \over
\sqrt{(x^2-a^2)(x^2+b^2)}}\;\left[(x+\sqrt{x^2-a^2})^{\nu}
+(x-\sqrt{x^2-a^2})^{\nu}\right]=2a^{\nu-\mu-1}b^{2\mu}\;\times
&\nonumber\\&
{\Gamma((\mu+\nu+1)/2) \Gamma((\mu-\nu+1)/2) \over \Gamma(\mu+1)}\;
_2F_1\left({\mu+\nu+1 \over 2},{\mu-\nu+1 \over 2};\mu+1;-{b^2
\over a^2}\right)
&\nonumber\\&
[{\rm Re}\;\mu>-1,\;{\rm Im}\;b=0,\;{\rm Re}\;(\mu+1)>
|{\rm Re}\;\nu|]
\label{beleven}&\end{eqnarray}

[No. 2.1.9.21 Ref.\ \cite{pru92a}, No. 2.16.21.1 Ref.\ \cite{pru86}],

\vspace{0.5cm}
\begin{eqnarray}&
\int_a^{\infty}dx\;\left[{(x+\sqrt{x^2-a^2})^{\nu} +(x-\sqrt{x^2-a^2})^{\nu}
\over \sqrt{x^2-a^2}}\right]x^{-\mu}K_{\mu}(bx)={\sqrt{\pi} \over
\Gamma(\mu+1/2)}\;\;\times
&\nonumber\\&
\left[2^{\nu-\mu}b^{\mu-\nu} \Gamma(\nu)\;\
B\left(\mu+{1 \over 2},{\nu \over 2}-\mu\right)\;_1F_2\left({1-\nu
\over 2};1-\nu,\mu-{\nu \over 2}+1;{a^2 b^2 \over 4}\right) + \right.
&\nonumber\\&
2^{-\nu-\mu-1} a^{2\nu}b^{\mu+\nu}\Gamma(-\nu)B\left(\mu+{1 \over 2},
-\mu-{\nu \over 2}\right)\;
_1F_2\left({1+\nu \over 2};1+\nu,\mu+{\nu \over 2}+1;{a^2 b^2 \over
4}\right)\;+
&\nonumber\\&
\left.2^{\mu-1}a^{\nu-2\mu}b^{-\mu}\Gamma\left(\mu+{\nu \over 2}\right)
\Gamma\left(\mu-{\nu \over 2}\right)\;_1F_2\left({1 \over 2}-\mu;
1-\mu-{\nu \over 2},1+{\nu \over 2}-\mu;{a^2 b^2 \over 4}\right)\right]
&\nonumber\\&
[{\rm Re}\;\mu>-1/2,\;b>0]
\label{bthirteen}&\end{eqnarray}

[No. 3.16.1.6 Ref.\ \cite{pru92a}, No. 2.16.3.7 Ref.\ \cite{pru86}],
\newline
and

\vspace{0.5cm}
\begin{eqnarray}&
\int_a^{\infty}dx\;\left[{(x+\sqrt{x^2-a^2})^{\nu}+(x-\sqrt{x^2-a^2})^{\nu}
\over \sqrt{x^2-a^2}}\right]x^{-\mu/2}\;K_{\mu}(b\sqrt{x})={2^{2\mu-2}
a^{\nu-\mu} \over b^{\mu}}\;\;\times
&\nonumber\\&
\Gamma\left({\mu+\nu \over 2}\right) \Gamma\left({\mu-\nu \over 2}\right)
\;_0F_3\left({1 \over 2},1-{\mu+\nu \over 2},1+{\nu-\mu \over 2};{a^2 b^4
\over 256}\right)-{2^{2\mu-4} a^{\nu-\mu+1} \over b^{\mu-2}}\;\;\times
&\nonumber\\&
\Gamma\left({\mu+\nu-1 \over 2}\right) \Gamma\left({\mu-\nu-1 \over 2}\right)
\;_0F_3\left({3 \over 2},{3-\mu-\nu \over 2},{3+\nu-\mu \over 2};{a^2 b^4
\over 256}\right)\;\;+
&\nonumber\\&
2^{-\mu-3\nu-1}a^{2\nu}b^{\mu+2\nu}\Gamma(-\nu) \Gamma(-\nu-\mu)\;
_0F_3\left(1+\nu,{1+\nu+\mu \over 2},1+{\nu+\mu \over 2};{a^2 b^4
\over 256}\right)\;\;+
&\nonumber\\&
2^{3\nu-\mu-1}b^{\mu-2\nu} \Gamma(\nu) \Gamma(\nu-\mu)\;
_0F_3\left(1-\nu,{1+\mu-\nu \over 2},1+{\nu-\mu \over 2};
{a^2 b^4 \over 256}\right)
&\nonumber\\&
[|{\rm arg}\;b|<\pi/4]
\label{bthirteena}&\end{eqnarray}

[No. 3.16.1.14 Ref.\ \cite{pru92a}, No. 2.16.8.9 Ref.\ \cite{pru86}].

Specifically, if we put $\mu=0$ in Eq.\ (\ref{bthirteen}), then we
find
\begin{eqnarray}&
\int_a^{\infty}dx\;{(x+\sqrt{x^2-a^2})^{\nu} +(x-\sqrt{x^2-a^2})^{\nu}
\over \sqrt{x^2-a^2}}\;K_0(b\sqrt{x})=a^{\nu}
K^2_{\nu/2}\left({ab \over 2}\right)\quad,
\label{bthirteenb}&\end{eqnarray}
or putting $\mu=1$, we get
\begin{eqnarray}&
\int_a^{\infty}dx\;{(x+\sqrt{x^2-a^2})^{\nu}+(x-\sqrt{x^2-a^2})^{\nu} 
\over \sqrt{x^2-a^2}}\;x^{-1}\;K_1(b\sqrt{x})=a^{\nu-1}\;\;\times
&\nonumber\\&
K_{(1+\nu)/2}\left({ab \over 2}\right)K_{\nu/2}\left({ab \over 2}\right)\;.
\label{bthirteenc}&\end{eqnarray}
If we put $\mu=0$ and $\nu=0$ in (177), then we obtain
\begin{eqnarray}&
\int_a^{\infty}dx\;{ K_0(b\sqrt{x}) \over \sqrt{x^2-a^2}}=
2[ker^2(b\sqrt{a/2})+kei^2(b\sqrt{a/2})]
&\nonumber\\&
\label{bthirteenf}&\end{eqnarray}

We can develop the general result given by Eqs. (\ref{bnine})
and (\ref{bten}) further by replacing $\nu$ by $n+\nu$ and multiplying
both equations by $(-1)^n \beta^n/n!$. Then by summing from $n=0$
to $\infty$ and applying the multiplication theorem for the MacDonald
function as given on p. 112 of Mangulis \cite{man65}, we get
\begin{eqnarray}&
\int_a^{\infty}dx\;{F(x) \over \sqrt{x^2-a^2}}\left[(x+\sqrt{x^2-a^2})^{\nu}
e^{-\beta (x+\sqrt{x^2-a^2})}+(x-\sqrt{x^2-a^2})^{\nu}
e^{-\beta (x-\sqrt{x^2-a^2})}\right]\;=
&\nonumber\\&
2a^{\nu}\int_0^{\infty}dt\;G(t)(1+2\beta/t)^{-\nu/2} K_{\nu}
(a\sqrt{t^2+2\beta t})\quad.
\label{bfourteen}&\end{eqnarray}
For the case of $\nu=0$, Eq.\ (\ref{bfourteen}) reduces to
\begin{eqnarray}&
\int_a^{\infty}dx\;{F(x) \over \sqrt{x^2-a^2}}\;e^{-\beta x}
\cosh(\beta \sqrt{x^2-a^2})=\int_0^{\infty}dt\;G(t)
K_0(a\sqrt{t^2+2\beta t})\quad.
\label{bfifteen}&\end{eqnarray}
In a similar manner, we can develop Eq.\ (\ref{bfive}) further.
This involves summing divergent integrals, the validity of which
has been discussed in Sec.\ 3. Thus, by multiplying both sides
of Eq.\ (\ref{bfive}) by $(-1)^n (\beta/2)^{2n+\nu}/n! \Gamma(n+\nu+1)$
and summing from $n=0$ to $\infty$, we eventually obtain after applying the
multiplication theorem:
\begin{eqnarray}&
\int_{2a}^{\infty}dx\;x^{\nu/2}(x-2a)^{\nu/2} J_{\nu}\left(\beta
\sqrt{x(x-2a)}\right)F(x)={a^{\nu+1/2} \beta^{\nu} \over 2^{\nu} \sqrt{\pi}}
\int_0^{\infty}dt\;e^{-at} G(t)\;\times
&\nonumber\\&
(t^2+\beta^2)^{-\nu/2-1/4} K_{\nu+1/2}\left(a\sqrt{t^2+\beta^2}\right)\;,
\label{bfifteena}&\end{eqnarray}
which is just a special case of Eq.\ (\ref{sevennine}) after
making the substitution $y=x^{1/2}(x-2a)^{1/2}$. This, therefore,
validates our technique of summing divergent integrals since Eq.\
(\ref{sevennine}) was obtained using convergent integrals. 

To complete this appendix we consider integrals of the form
\begin{eqnarray}&
I=\int_0^{\infty}dx\;erfc(ax^{-1/2})F(x) \quad,
\label{bsixteen}&\end{eqnarray}
which can be transformed into
\begin{eqnarray}&
I=\int_0^{\infty}dy\;y^{-1}e^{-ay}G(y^2)\quad.
\label{bseventeen}&\end{eqnarray}
An example arising from the above results is 
\begin{eqnarray}&
\int_0^{\infty}dx\;{\rm erfc}(ax^{-1/2})\; (x^2+p^2)^{\nu}=
{\sqrt{\pi}\;p^{2\nu+1}\Gamma(-\nu-1/2) \over 4 \Gamma(-\nu)}-
{\sqrt{\pi}\;p^{2\nu+1/2}a \over 2^{3/2} \Gamma(3/4)\Gamma(-\nu)}\;\times
&\nonumber\\&
\Gamma\left(-\nu-{1 \over 4}\right)\;_2F_3\left({1 \over 4},
-{\nu+1 \over 4};{1 \over 2},{3 \over 4},{5 \over 4};-{a^4 \over 64 p^2}\right)
+{p^{2\nu} a^2 \over 4}
\;_2F_3\left({1 \over 2},-\nu;{3 \over 4},{5 \over 4},{3 \over 2};-
{a^4 \over 64 p^2}\right)-
&\nonumber\\&
\sqrt{{\pi \over 2}}\;{p^{2\nu-1/2}a^3 \over 6}\;
{\Gamma(1/4-\nu) \over \Gamma(1/4) \Gamma(-\nu)}
\;_2F_3\left({3 \over 4},{1 \over 4}-\nu;{5 \over 4},
{3 \over 2},{7 \over 4};-{a^4 \over 64 p^2}\right)
&\nonumber\\&
[{\rm Re}\;\nu<-1/2,\;{\rm Im}\;p=0]
\label{bseventeena}&\end{eqnarray}

[No. 2.1.5.1 Ref.\ \cite{pru86}, No. 2.12.9.4 Ref.\ \cite{pru92a}].

\section{Appendix C}
Here we establish the Feynman integral from its more
general counterpart given as No. 2.2.6.1 in Prudnikov et al
\cite{pru86a} and then develop further integral results from it.
This integral, which appears so prominently
in quantum electrodynamics, can be written generally \cite{pok89}
as
\begin{equation}
\int_0^1dx\;{x^{\alpha-1}\;(1-x)^{\beta-1} \over
(ax+b(1-x))^{\alpha+\beta}}={\Gamma(\alpha) \Gamma(\beta)
\over \Gamma(\alpha+\beta) a^{\alpha}\;b^{\beta}} \quad,
\label{cone}\end{equation}
where Re $\alpha >0$ and Re $\beta >0$. These conditions will
apply throughout this appendix. To establish the above result,
we shall show that
\begin{eqnarray}&
I=\int_a^bdx\;{(x-a)^{\alpha-1}\;(b-x)^{\beta-1} \over
(cx+d)^{\gamma}}=(b-a)^{\alpha+\beta-1}(ac+d)^{\gamma}\;
B(\alpha,\beta)\;\;\times
&\nonumber\\&
\;_2F_1\left(\alpha,\gamma;\alpha+\beta;{c(a-b) \over ac+d}\right) \;,
\label{ctwo}&\end{eqnarray}
where Re $\alpha$,  $\beta > 0$ and $|{\rm arg} ((d+cb)/(d+ca))| < \pi$.
Eq.\ (\ref{ctwo}) is simply No.\ 2.2.6.1 from Ref.\ \cite{pru86a} with
$\gamma$ replaced by $-\gamma$. We shall primarily be interested in the
case of Re $\gamma > 0$ even though the case of Re $\gamma < 0$ can
be established by analytic continuation.

With the aid of the integral representation of the gamma function
Eq.\ (\ref{ctwo}) can be written as
\begin{equation}
I= {(b-a)^{\alpha+\beta-1} \over \Gamma(\gamma)}\;\int_0^1dx\;
x^{\alpha-1}\;(1-x)^{\beta-1}
\int_0^{\infty}dt\;t^{\gamma-1}\;e^{-(fx+g)t} \quad,
\label{cthree}\end{equation}
where $f=c(b-a)$ and $g=ac+d$. Expanding the part of the exponential
with $fxt$ into a series and then interchanging the order of the
integrations and the sum, we find that Eq.\ (\ref{cthree}) yields
\begin{eqnarray}&
I={(b-a)^{\alpha+\beta-1} \over \Gamma(\gamma)}\;B(\alpha,\beta)\;
\sum_{n=0}^{\infty} {(-1)^nf^n \over n!} {\Gamma(n+\alpha)
\over \Gamma(\alpha)}{\Gamma(\alpha+\beta) \over \Gamma(n+\alpha+\beta)}\;
{\Gamma(n+\gamma) \over g^{\gamma}}\;\;=
&\nonumber\\&
(b-a)^{\alpha+\beta-1}\;B(\alpha,\beta)\;
_2F_1(\alpha,\gamma;\alpha+\beta;-f/g) \quad,
\label{cfour}&\end{eqnarray}
which is just Eq.\ (\ref{ctwo}). To obtain the Feynman integral
put $a=0$, $b=1$, $c=a-b$, $d=b$ and $\gamma=\alpha+\beta$. Then
the $_2F_1$ hypergeometric function reduces to a
$_1F_0$ hypergeometric function, which in turn yields an algebraic
function as given on p. 453 of Ref.\ \cite{pru90}. Thus one finally
arrives at Eq.\ (\ref{cone}).

We shall consider the more general form of Eq.\ (\ref{cone}) with
$\alpha+\beta$ replaced by $\gamma$ in the power of the
denominator shortly. For now if we put $\alpha$ equal $n+\alpha$
with $n$ a non-negative integer and multiply both sides by 
$s^n/n!$, then we obtain after summing from $n=0$ to $\infty$
and interchanging the order of the integration and summation
\begin{equation}
\int_0^{1}dx\;{x^{\alpha-1}\;(1-x)^{\beta-1} \over
(ax+b(1-x))^{\alpha+\beta}}\;e^{sx/(ax+b(1-x))}\;=
B(\alpha,\beta)\;a^{-\alpha}\;b^{-\beta}\;
_1F_1(\alpha;\alpha+\beta;s/a) \quad.
\label{cfive}\end{equation}
On the other hand if we had multiplied both sides by $(-1)^ns^n/n!
\Gamma(n+\gamma+1)$, then we would obtain
\begin{eqnarray}&
\int_0^1dx\;{x^{\alpha-\gamma/2-1}\;(1-x)^{\beta-1} \over
(ax+b(1-x))^{\alpha+\beta-\gamma/2}}\;
J_{\gamma}\left(2\sqrt{{sx \over ax+b(1-x)}}\right)=
{B(\alpha,\beta) \over \Gamma(\gamma+1)}\;\;\times
&\nonumber\\&
{s^{\gamma/2}\over b^{\beta}\;a^{\alpha}}\;_1F_2\left(\alpha;\gamma+1,
\alpha+\beta;-{s \over a}\right) \quad.
\label{csix}&\end{eqnarray}
Alternatively, dropping the phase factor of $(-1)^n$ in the process
yields
\begin{eqnarray}&
\int_0^1dx\;{x^{\alpha-\gamma/2-1}\;(1-x)^{\beta-1} \over
(ax+b(1-x))^{\alpha+\beta-\gamma/2}}\;
I_{\gamma}\left(2\sqrt{{sx \over ax+b(1-x)}}\right)=
{B(\alpha,\beta) \over \Gamma(\gamma+1)}\;\;\times
&\nonumber\\&
{s^{\gamma/2}\over b^{\beta}\;a^{\alpha}}\;_1F_2\left(\alpha;\gamma+1,
\alpha+\beta;{s \over a}\right) \quad.
\label{cseven}&\end{eqnarray}

If we put $\gamma=\alpha$ in Eq. \ref{csix}, then by utilising 
tabulated results for $_1F_2$ hypergeometric functions on p. 608 of
Ref.\ \cite{pru90} we get
\begin{eqnarray}&
\int_0^1dx\;{x^{\alpha/2-1}\;(1-x)^{\beta-1} \over 
(ax+b(1-x))^{\alpha/2+\beta}}\;J_{\alpha}\left(2\sqrt{{sx \over
ax+b(1-x)}}\right)={2^{1+\beta-\alpha}\;s^{(1-\alpha)/2} \Gamma(\beta) 
\over a^{1/2}\;b^{\beta}}\;\;\times
&\nonumber\\&
\left(2 \alpha J_{\alpha+\beta-1}\left(2\sqrt{{s \over a}}
\right)\;s_{\alpha-\beta-1,\alpha+\beta-2}\left(2\sqrt{{s \over a}}
\right)-J_{\alpha+\beta-2}\left(2\sqrt{{s \over a}}\right)\;
s_{\alpha-\beta,\alpha+\beta-1}\left(2\sqrt{{s \over a}}\right)\right)
\quad,
\label{ceight}&\end{eqnarray}
and for the case where $\alpha+\beta=3/2$, the above equation yields
\begin{eqnarray}&
\int_0^1dx\;{x^{\alpha/2-1}\;(1-x)^{1/2-\alpha} \over
(ax+b(1-x))^{(3\alpha-1)/2}}\;J_{\alpha}\left(2\sqrt{{sx \over
ax+b(1-x)}}\right)={2^{2-2\alpha}\;\Gamma(3/2-\alpha) \over
b^{3/2-\alpha}\;s^{\alpha/2} \sqrt{\pi}}
&\nonumber\\&
\left[\Gamma(2\alpha-1)\cos(\pi \alpha)+ S\left(2\sqrt{{s \over a}},
2\alpha-1\right)\right] \quad.
\label{cnine}&\end{eqnarray}
In Eq.\ (\ref{ceight}), $s_{\mu,\nu}(z)$ represents the Lommel function
while in Eq.\ (\ref{cnine}), $S(z,\nu)$ represents the generalised Fresnel
sine integral and Re $\alpha < 1/2$.

If we put $\gamma=\alpha-2$ in Eq.\ ({\ref{cseven}), then we find using
the results on p. 595 of Ref.\ \cite{pru90} that
\begin{eqnarray}&
\int_0^1dx\;{x^{\alpha/2}\;(1-x)^{\beta-1} \over (ax+b(1-x))^{\beta+
\alpha/2+1}}\;I_{\alpha-2}\left(2\sqrt{{sx \over ax+b(1-x)}}\right)=
{(\alpha-1)a^{(\beta-\alpha-1)/2}\;\Gamma(\beta) \over b^{\beta}
s^{(\beta +1)/2}}\;\;\times
&\nonumber\\&
\left[I_{\alpha+\beta-1}\left(2\sqrt{{s \over a}}\right) + {(s/a)^{1/2}
\over 2(\alpha-1)}\;I_{\alpha+\beta}\left(2 \sqrt{{s \over a}}\right)
\right] \quad,
\label{cten}&\end{eqnarray}
or if we put $\gamma=\alpha-1/2$ and $\beta=\alpha$, then we find
\begin{eqnarray}&
\int_0^1dx\;{x^{\alpha/2-3/4}\;(1-x)^{\alpha-1} \over
(ax+b(1-x))^{3\alpha /2+1/4}}\;I_{\alpha-1/2}\left(2 \sqrt{sx \over
ax+b(1-x)}\right)= {2^{2\alpha-1}\;a^{-1/2} \over b^{\alpha}\;
s^{\alpha/2-1/4}} \;\;\times
&\nonumber\\&
\Gamma(\alpha+1/2)\;B(\alpha,\alpha)\;I^2_{\alpha-1/2}\left(
\sqrt{{s \over a}}\right) \quad.
\label{coneone}&\end{eqnarray}
On the other hand, putting $\alpha=1/2$, $\gamma=1-c$ and $\beta=
c-1/2$ yields
\begin{eqnarray}&
\int_0^1dx\;{x^{c/2-1}\;(1-x)^{c-3/2} \over (ax+b(1-x))^{3c/2-1/2}}\;
I_{1-c}\left(2 \sqrt{{sx \over ax+b(1-x)}}\right)=
{\pi B(1/2,c-1/2) \over \Gamma(2-c) \sin(\pi c)}\;\;\times
&\nonumber\\&
{(1-c)\;s^{(1-c)/2} \over b^{c-1/2}\;a^{1/2}}\;
I_{1-c}\left(\sqrt{{s \over a}}\right)\; I_{c-1}\left(\sqrt{{s \over a}}
\right) \quad,
\label{conetwo}&\end{eqnarray}
the latter result being valid for $0<$ Re $c < 2$.

We should also note that particular results for Eqs. (\ref{csix}) and 
(\ref{cseven}) can be obtained by utilising the tabulated results of
$_1F_2$ hypergeometric functions in Prudnikov et al \cite{pru90}
with specific indices. For example, we find
\begin{eqnarray}&
\int_0^1dx\; {x^{-7/8}\;(1-x)^{-3/4} \over (ax+b(1-x))^{3/8}}\;
J_{1/4}\left(2 \sqrt{{sx \over ax+b(1-x)}}\right)=
{2\Gamma(1/4)\; a^{1/4} \over b^{1/4} s^{3/8}}\;
C\left(2\sqrt{{s \over a}}\right) \quad,
\label{conethree}&\end{eqnarray}
\begin{eqnarray}&
\int_0^1dx\;{x^{-1} \over ax+b(1-x)}\;\sin\left(2\sqrt{{sx \over
ax+b(1-x)}}\right)={2 \over b}\;
Si\left(2\sqrt{{s \over a}}\right) \quad,
\label{conefour}&\end{eqnarray}
\begin{eqnarray}&
\int_0^1dx\;{x^{-1}\;(1-x)^{-1/2} \over (ax+b(1-x))^{1/2}}\;
I_1\left(2\sqrt{{sx \over ax+b(1-x)}}\right) ={\pi s^{1/2}
\over \sqrt{a b}}\;\left(I_0^2\left(\sqrt{s \over a}\right)-
I_1^2\left(\sqrt{{s \over a}}\right)\right)\;,
\label{conefive}&\end{eqnarray}
and 
\begin{eqnarray}&
\int_0^1dx\;{x^{-1/4} \over (ax+b(1-x))^{7/4}}\;I_{1/2}
\left(2\sqrt{{sx \over ax +b(1-x)}}\right)=
{1 \over \sqrt{\pi}\;b s^{3/4}}\;\left[
\cosh\left(2\sqrt{{s \over a}}\right)-1\right] \quad.
\label{conesix}&\end{eqnarray}
We add that many more results can be obtained by following the
above procedure.

Now returning to the generalised Feynman integral, i.e. Eq.
(\ref{cone}), we note that if we put $\alpha=n+\alpha+1$,
$\beta=n+\beta$ and multiply both sides by $c^{2n+\alpha}/
2^{2n+\alpha}n!\Gamma(n+\alpha+1)$, then we obtain the following
result after interchanging the order of the sum from 0 to
$\infty$ and integration
\begin{eqnarray}&
\int_0^1dx\;{x^{\alpha/2}(1-x)^{\beta-\alpha/2-1} \over 
(ax+b(1-x))^{\beta+1}}\;I_{\alpha}\left({c \sqrt{x(1-x)} \over
ax+b(1-x)}\right)={c^{\alpha}\; \Gamma(\beta) \over 2^{\alpha}
a^{\alpha+1}\;b^{\beta}\;\Gamma(\alpha+\beta+1)}\;\;\times
&\nonumber\\&
_1F_2\left(\beta;{\alpha+\beta+1 \over 2},{\alpha+\beta \over 2}+1;{c^2 \over 16ab}
\right) \quad,
\label{coneseven}&\end{eqnarray}
while if we multiply by $(-1)^n c^{2n+\alpha}/2^{2n+\alpha}n!
\Gamma(n+\alpha+1)$ instead, we get
\begin{eqnarray}&
\int_0^1dx\;{x^{\alpha/2}\;(1-x)^{\beta-\alpha/2-1} \over
(ax+b(1-x))^{\beta+1}}\;J_{\alpha}\left({c \sqrt{x(1-x)} \over
ax+b(1-x)}\right)={c^{\alpha}\;\Gamma(\beta) \over 2^{\alpha}
a^{\alpha+1}\;b^{\beta}\;\Gamma(\alpha+\beta+1)}\;\;\times
&\nonumber\\&
_1F_2\left(\beta;{\alpha+\beta+1 \over 2},{\alpha+\beta \over 2}+1;
-{c^2 \over 16ab}\right) \quad.
\label{coneeight}&\end{eqnarray}
Some results emanating from the above are
\begin{eqnarray}&
\int_0^1dx\;{x^{(\gamma-3)/2}\;(1-x)^{(\gamma+1)/2} \over
(ax+b(1-x))^{\gamma+1}}\;I_{\gamma-3}\left({c \sqrt{x(1-x)}
\over ax+b(1-x)}\right)={2^{\gamma}\;\Gamma(\gamma) \Gamma(\gamma-1/2)
\over a^{\gamma/2-5/4}\;b^{\gamma/2+3/4}\;\Gamma(2\gamma-2)}\;\times
&\nonumber\\&
c^{-3/2}\left[I_{\gamma-3/2}\left({c \over 2\sqrt{ab}}\right)- {c \over 8(\gamma-1)
\sqrt{ab}}\;I_{\gamma-1/2}\left({c \over 2\sqrt{ab}}\right)\right] \quad,
\label{conenine}&\end{eqnarray}
\begin{eqnarray}&
\int_0^1dx\;{x^{1/2}\;(1-x)^{-1/2} \over (ax+b(1-x))^2}\;
I_1\left({c\sqrt{x(1-x)} \over ax+b(1-x)}\right)=
{2 \over a c}\left[\cosh\left({c \over 2\sqrt{ab}}\right)-1\right] \quad,
\label{ctwozero}&\end{eqnarray}
and
\begin{eqnarray}&
\int_0^1dx\;{x^{1/2}\;(1-x)^{-1/2} \over (ax+b(1-x))^2}\;
J_1\left({c \sqrt{x(1-x)} \over ax+b(1-x)}\right)=
\sqrt{{\pi \over c}}\; a^{-5/4}b^{-1/4}{\bf H}_{1/2}\left(
{c \over 2\sqrt{ab}}\right) \quad .
\label{ctwoone}&\end{eqnarray}
In Eq.\ (\ref{ctwoone}), ${\bf H}_{\nu}(z)$ refers to a Struve
function. In addition, by putting $\alpha=\beta-1$ in Eq.\ (\ref{coneeight}),
one obtains
\begin{eqnarray}&
\int_0^1dx\;{x^{(\beta-1)/2}(1-x)^{(\beta-1)/2} \over (ax+b(1-x))^{\beta+1}}\;
J_{\beta-1}\left({c\sqrt{x(1-x)} \over ax+b(1-x)}\right)={2^{1-\beta}c^{-1/2}
\sqrt{\pi} \over (ab)^{\beta/2+1/4}}\;J_{\beta-1/2}\left({c \over
2\sqrt{ab}}\right).
\label{ctwozeroa}&\end{eqnarray}

By utilising similar techniques to those described in this
appendix, we can develop further results based on the Feynman
integral such as
\begin{eqnarray}&
\int_0^1 dx\;{x^{\alpha-\gamma/2-1}\;(1-x)^{\beta-\gamma/2-1}
\over (ax+b(1-x))^{\alpha+\beta-\gamma}}\;I_{\gamma}
\left({2 \sqrt{cx(1-x)} \over ax+b(1-x)}\right)={\Gamma(\alpha)
\Gamma(\beta) c^{\gamma/2} \over \Gamma(\gamma+1) 
\Gamma(\alpha+\beta) a^{\alpha} \;b^{\beta}} \;\;\times
&\nonumber\\&
_2F_3\left(\alpha,\beta;{\alpha+\beta \over 2},
{\alpha+\beta+1 \over 2}, \gamma;{c \over 4ab}\right) \quad,
\label{ctwotwo}&\end{eqnarray}
\begin{eqnarray}&
\int_0^1 dx\;{x^{\alpha-1}\;(1-x)^{-\alpha} \over
(ax+b(1-x))}\;_0F_3\left(\alpha,1-\alpha, {1 \over 2};
{-z^4 x(1-x) \over 16(ax+b(1-x))^2}\right)=
&\nonumber\\&
{\Gamma(\alpha) \Gamma(1-\alpha) \over a^{\alpha}
b^{1-\alpha}}\;ber\left({z \over (ab)^{1/4}}\right)\quad, \quad {\rm Re}
\;\alpha<1 \quad .
\label{ctwothree}&\end{eqnarray}
A similar result to Eq.\ (\ref{ctwothree}) can be found for  the Kelvin function,
$bei(z)$, by using No. 8.564 from Gradshteyn and Ryzhik \cite{gra80}.
Finally, Eq.\ (\ref{ctwotwo}) yields
interesting results of its own, two of which are
\begin{eqnarray}&
\int_0^1dx\;x^{(\alpha-\beta)/2-1}\;(1-x)^{(\beta-\alpha)/2-1}\;
I_{\alpha+\beta}\left({2 \sqrt{cx(1-x)} \over ax+b(1-x)}\right)=
{B(\alpha,\beta)\;c^{(\alpha+\beta)/2} \over \Gamma(\alpha+\beta+1)}\;\;\times
&\nonumber\\&
a^{-\alpha}\;b^{-\beta}\; _1F_1\left(\alpha;\alpha+\beta;\sqrt{ c \over ab}\right)
\;_1F_1\left(\alpha;\alpha+\beta;-\sqrt{c \over ab}\right) \quad,
\label{ctwofive}&\end{eqnarray}
and 
\begin{eqnarray}&
\int_0^1dx\;{x^{-1}(1-x)^{-1/2} \over (ax+b(1-x))^{1/2}}\;
I_{2\alpha}\left({2\sqrt{cx(1-x)} \over ax+b(1-x)}\right)=
{2^{4\alpha-2} c^{1/2} \over a^{1/2}\;b}\;\;B(\alpha+1/4,\alpha+3/4)\;\times
&\nonumber\\&
B(\alpha,\alpha+1/2)\;I_{\alpha-1/4}\left( \sqrt{{c \over 4ab}}\right)
I_{\alpha-3/4}\left(\sqrt{{c \over 4ab}}\right) \quad.
\label{ctwosix}&\end{eqnarray}

\end{document}